\pdfoutput=1
\documentclass[12pt]{iopart}
\usepackage{url,doi}

\expandafter\let\csname equation*\endcsname\relax
\expandafter\let\csname endequation*\endcsname\relax
\usepackage{amsmath, amsthm, amsfonts}

\usepackage{epic}
\usepackage{graphicx}
\usepackage[dvips]{epsfig}
\usepackage[usenames, dvipsnames]{color}
\usepackage{enumitem}
\usepackage{multicol}
\usepackage{tabularx}
\usepackage{paracol}
\usepackage{caption}
\usepackage{csquotes}
\usepackage{xcolor}

\theoremstyle{plain} 
  \newtheorem{theo}{Theorem}[section]

  \newtheorem{coro}[theo] {Corollary}
  \newtheorem{prop}[theo] {Proposition}

\theoremstyle{definition}
  
  \newtheorem{remark}[theo] {Remark}
  
  \newtheorem{assumption}[theo] {Assumption}

\newcommand{\N}{{\mathbb{N}}}

\newcommand{\R}{{\mathbb{R}}}  

\newcommand{\F}{{\mathbb{F}}}
\newcommand{\Y}{{\mathbb{Y}}}

\newcommand{\bF}{{\mathbf{F}}}
\newcommand{\bU}{{\mathbf{U}}}
\newcommand{\bG}{{\mathbf{G}}}

\newcommand{\calD}{{\cal D}}
\newcommand{\calB}{{\cal B}}

\newcommand{\calN}{{\cal N}}
\newcommand{\calX}{{\cal X}}
\newcommand{\calY}{{\cal Y}}
\newcommand{\calU}{{\cal U}}
\newcommand{\calW}{{\cal W}}

\def\paren#1{\left( #1 \right)}

\def\dupair#1{\left\langle #1 \right\rangle}

\def\inprod#1{\langle #1 \rangle}

\newcolumntype{L}[1]{>{\raggedright\arraybackslash}p{#1}}
\newcolumntype{C}[1]{>{\centering\arraybackslash}p{#1}}
\newcolumntype{R}[1]{>{\raggedleft\arraybackslash}p{#1}}

\hypersetup{
  colorlinks,
  citecolor=Blue,
  linkcolor=Blue,
  urlcolor=Blue}

\colorlet{blue}{black}  	

\begin{document}

\title[Landweber-Kaczmarz in time-dependent inverse problems]{Landweber-Kaczmarz for parameter identification \\ in time-dependent inverse problems:\\All-at-once versus Reduced version}
\author{\bf Tram Thi Ngoc Nguyen}
\address{Alpen-Adria-Universit\"at Klagenfurt, Universit\"atstra{\ss}e 65-67, A-9020 Klagenfurt, Austria}
\ead{Tram.Nguyen@aau.at}

\date{\today}


\begin{abstract}
In this study, we consider a general time-space system, \textcolor{blue}{whose model operator and observation operator are locally Lipschitz continuous}, over a finite time horizon and parameter identification by using Landweber-Kaczmarz regularization. The problem is investigated in two different modeling settings: An All-at-once and a Reduced version, together with two observation scenarios: continuous and discrete observations. Segmenting the time line into several subintervals leads to the idea of applying the Kaczmarz method. A loping strategy is incorporated into the method to yield the loping Landweber-Kaczmarz iteration.

\end{abstract}

\section{Introduction} \label{secintro}
As a motivating prototype example, we consider parameter identification from the following system
\begin{align}
& \dot{u}=\Delta u - u^3 + \theta \qquad\hphantom{.} (t,x) \in (0,T)\times \Omega \label{model-intro1}\\
& u(0)=u_0 \qquad\qquad\quad \text{ } x \in \Omega\label{model-intro2} \\
& y=Cu \qquad\qquad\qquad (t,x) \in (0,T)\times \Omega\label{model-intro3}
\end{align}
with $\Omega \subset \R^d$. $\dot{u}$ denotes the first order time derivative of $u$ and the right hand side includes the nonlinearity $\Phi(u):=u^3$.
This equation is equipped with the initial condition (\ref{model-intro2}) and possibly further Dirichlet or Neumann boundary conditions on $(0,T)\times \partial\Omega$. In (\ref{model-intro3}), the measured data $y$ is obtained from a linear observation, this means $C$ is a linear operator. In this evolution system, $u$ and $\theta$ are two unknowns.

Parabolic PDEs with cubic power nonlinearity arise in many applications, and we selectively mention some falling into this category
\begin{itemize}
\item $\Phi(u)=u(1-u^2)$: 
Ginzburg-Landau equations of superconductivity \cite{BronsardStoh}, Allen-Cahn equation for phase separation process in a binary metallic alloy \cite{AllenCahn, Nepomnyashchy}, Newell-Whitehead equation for convection of fluid heated from below \cite{GildingKersner}.

\item $\Phi(u)=u^2(1-u)$:
Zel'dovich equation in combustion theory \cite{GildingKersner}.

\item $\Phi(u)=u(1-u)(u-\alpha), 0<\alpha<1$: Fisher's model for population genetics \cite{Pao}, Nagumo equation for bistable transmission lines in electric circuit theory \cite{NagumoYoshizawaArinomoto}. 
\end{itemize}

To the best of our knowledge, the state-of-the-art research on this type of problem in Sobolev space framework is limited to the power $\gamma\leq 1+\frac{4}{d}, d\geq 3$, where $d$ is the space dimension \cite{Kaltenbacher}. The reason for this constraint is that, the growth conditions used in proving well-definedness and differentiability of the forward operators prevent higher nonlinearity. Therefore, in this study, staying in the Sobolev space framework, we aim at increasing the nonlinearity to the power of $3$, or even higher, in order to be compatible with those applications. For this purpose, we plan to construct appropriate function spaces for the forward operators and impose on them relevant assumptions. The main assumption we rely on is basically the local Lipschitz continuity condition. Recently, some authors \cite{TuanTrong, HaoDucThang} also employ local Lipschitz continuity as the key ingredient for the research on backward parabolic problems, but rather in a semigroup framework than in the Sobolev space framework, e.g., of \cite{Roubicek}, which we rely on to get even somewhat weaker conditions. \textcolor{blue}{Beyond increasing nonlinearity as compared to \cite{Kaltenbacher}, this paper provides explicit formulas for Hilbert space adjoints that makes our method more computationally efficient, and G\^ ateaux as well as Fr\' echet derivatives obtained in this work yield more knowledge on the properties of the nonlinear model and on the convergence performance of iterative methods.}

The fundamental idea of our approach is Landweber regularization \cite{HankeNeubauerScherzer}. Let $F$ be a nonlinear operator mapping between Hilbert spaces $\calX$ and $\calY$, and that has continuous G\^ateaux derivative. Let the problem under consideration be formulated by the model $F(x)=y$ and in case of noisy data, only $y^\delta$ as an approximation to $y$ is available. Starting at an initial point $x_0$, the Landweber iteration is defined by
\begin{align}
 x_{k+1}= x_k - \mu_k F'(x_k)^\star(F(x_k)-y^\delta) \qquad k\in\N_0,
\end{align}
where $F'(x)$ is the derivative of $F$ at $x$ and $F'(x)^\star$ is its adjoint. Note that $F$ is not necessarily Fr\'echet differentiable, just needs to be G\^ateaux differentiable with bounded linear derivative.
If $F=(F_0,\ldots, F_{n-1}):\bigcap\limits^{n-1}_{i=0} \calD(F_i) \subset \calX \rightarrow \calY^n$ and $y^\delta=(y^\delta_0,\ldots,y^\delta_{n-1})$, the Landweber-Kaczmarz method \cite{HaltmeierKowarLeitaoScherzer, HaltmeierLeitaoScherzer, KowarScherzer} reads as follows
\begin{align} \label{KaczmarzIntro}
x_{k+1}= x_k -\mu_k F'_{j(k)}(x_k)^\star(F_{j(k)}(x_k)-y_{j(k)}^\delta)\qquad k\in\N_0
\end{align}
with $j(k)=k-n\lfloor k/n\rfloor$, the largest integer lower or equal to $k$. It is recognizable that the Landweber-Kaczmarz method applies the Landweber iterations in a cyclic manner. Basing on this, we propose the idea of segmenting the time horizon $(0,T)$ into several subintervals, building up several corresponding forward operators and applying the Landweber-Kaczmarz method. 

Being more realistic than the model (\ref{model-intro1})-(\ref{model-intro3}) with observation on all of $(0,T)$, in practice the unknown parameter is recovered from experimental techniques that in some cases, limit the measurement only to some particular time points, see, e.g.,  \cite{Pedretscher, Boiger} for just two out of many examples from material science and system biology, respectively. Therefore, beside the continuous observation, we also desire to cover the discrete observation case in this study.

The paper is outlined as follows. Section \ref{model} introduces the general model for the state-space problem and discusses the function space setting as well as the reconstruction method being used. In the next two sections, we formulate the problem in two different settings: An All-at-once and a Reduced version. In each version, we take into account both cases of continuous and discrete observations. \textcolor{blue}{Section \ref{secdiscuss} compares the two modeling settings and discusses time-dependent parameter identification}. Section \ref{secAlgorithm} is dedicated to deriving the algorithm and collecting some convergence results for the Landweber-Kaczmarz method. In Section \ref{secEx}, we investigate an example by initially verifying the conditions proposed in Sections \ref{secAAO} and \ref{secReduced}, then running some numerical experiments. Finally, Section \ref{secOutlook} concludes the work and sketches some ideas for further potential research.


\section{Mathematical model} \label{model}
We consider the following state-space system
\begin{align}
& \dot{u}(t) = f(t,u(t),\theta) \quad t\in (0,T) \qquad\qquad\qquad\qquad \label{model-1}\\ 
& u(0)=u_0(\theta) \label{model-2}\\
& y(t) = g(t,u(t),\theta)  \quad t\in (0,T) \label{model-3}\\
&\text{or} \nonumber\\
& y_i = g_i(u(t_i),\theta)  \qquad i=1\ldots n \label{model-4}
\end{align}
on $\Omega \subset \R^d$, where $f$
is a nonlinear function
 and additional observation data $y$ or $y_i$ are obtained from continuous or discrete measurement as in (\ref{model-3}) or (\ref{model-4}), respectively. In the general case, $g, g_i$ may be linear or nonlinear. \textcolor{blue}{In particular, observations may be partial only, such as boundary traces on $\partial\Omega$ or a part of it.}

The model operator and observation operators map between the function spaces
\begin{align}
& f: (0,T)\times V\times \calX \rightarrow V^*\\
& g: (0,T)\times V\times \calX \rightarrow Z \quad\text{or}\quad g_i: V\times \calX \rightarrow Z, \quad
\end{align}
where $\calX, Z, H$ and $V$ are separable Hilbert spaces and $V \hookrightarrow H \hookrightarrow V^*$ form a Gelfand triple.
Moreover, we assume that $f, g$ meet the Caratheodory mapping condition.

The initial condition is supposed to map to the sufficiently smooth image space
\begin{align}
&\hphantom{hh} u_0:\calX \rightarrow V \qquad\qquad\qquad\qquad\qquad\qquad\qquad\qquad
\end{align}
as a condition to attain some regularity results for the solution to the problem (\ref{model-1})-(\ref{model-2}).

Between $V$ and $V^*$, the Riesz isomorphism
\begin{align*}
& I:V^*\rightarrow V,\quad \dupair{u^*,Iv^*}_{V^*,V}=(u^*,v^*)_{V^*} \qquad\qquad
\end{align*}
and
\begin{align*}
& \tilde{I}:V^*\rightarrow V,\quad (\tilde{I}u^*,v)_V=\dupair{u^*,v}_{V^*,V} \qquad\qquad\quad
\end{align*}
are used to derive the adjoints. $\tilde{I}$ as defined above exists as one can choose, for example, $\tilde{I}=D^{-1}$, where $D$ is the Riesz isomorphism
\begin{align*}
& D:V\rightarrow V^*,\quad \dupair{Du,v}_{V^*,V}=(u,v)_V. \qquad\qquad\hphantom{hh}
\end{align*}
Here, $(.,.)$ and $\inprod{.,.}$ with the subscripts indicate the inner products and the dual parings, respectively. The notations $D$ and $I$ refer to the spatial differential and integration operators in the context of parabolic differential equations. We also distinguish the superscript $^*$, the Banach space adjoint, and $^\star$, the Hilbert space adjoint, which is an ingredient for the iterative methods considered here. 

For fixed $\theta$, $f$ and $g$ as defined above induce Nemytskii operators \cite[Section 4.3]{Troltzsch} on the function space $\calU$. This function space will be expressed later according to which problem setting is being dealt with, i.e., All-at-once or Reduced setting. However, they map into the same image space $\calW$ and observation space $\calY$
\begin{equation}
\calW=L^2(0,T;V^*), \qquad \calY=L^2(0,T;Z), \qquad\ \label{spaceYW}
\end{equation} 
which are Hilbert spaces.
Therefore, we can investigate the problem in the Hilbert space framework, provided that the corresponding argument spaces of the forward operators in the two settings are Hilbert spaces as well. 

In reality, we do not have access to the exact data. The experimental data always contains some noise of which we assume to know just the noise level. The noise perturbing the system is present both on the right hand side of the model equation and in the observation, and is denoted respectively by $w^\delta$ and $z^\delta$. When the measurement is a collection of data at discrete observation time points, the corresponding noise added to each observation is $z_i^\delta$. Altogether, we can formulate the noisy system
\begin{align}
& \dot{u}(t) = f(t,u(t),\theta)+w^\delta(t) \quad t\in (0,T) \qquad\qquad\quad \\ 
& u(0)=u_0(\theta) \\
& y^\delta(t) = g(t,u(t),\theta)+z^\delta(t)  \quad t\in (0,T) \\
&\text{or} \nonumber\\
& y^\delta_i = g_i(u(t_i),\theta)+z_i^\delta  \quad\qquad\hphantom{l} i=1\ldots n. 
\end{align}
Correspondingly, those additive noise terms live in the function spaces
\begin{equation}
w^\delta \in \calW, \qquad z^\delta \in\calY, \qquad z_i^\delta \in Z \qquad i=1\ldots n \qquad\qquad
\end{equation} 
and are supposed to satisfy
\begin{equation*}
\|w^\delta\|_\calW \leq \delta_w, \qquad \|z^\delta\|_\calY \leq \delta_z, \qquad \|z_i^\delta\|_Z \leq \delta_i \quad\hphantom{h} i=1\ldots n \qquad
\end{equation*}
with the noise levels $\delta_w,\delta_z,\delta_i > 0$.

This paper does not attempt to research uniqueness of the exact solution $(u^\dagger,\theta^\dagger)$. Nevertheless, we refer to the book \cite{BarbaraNeubauerScherzer}, which exposes some general results for this important question based on the tangential cone condition together with the assumption of a trivial null space of $F'(x^\dagger)$ in some neighborhood of $x^\dagger$. To verify the later condition in some concrete time-space problems, one can find detailed discussions, e.g., in the book by Isakov \cite{Isakov}.

\section{All-at-once setting} \label{secAAO}
In this section, we recast the system into a form which allows solving for both state $u$ and parameter $\theta$ simultaneously. For this purpose, we define the forward operator
\begin{align} \label{AAO-model}
&\F: \calU \times \calX  \rightarrow \calW \times V \times \calY,  \qquad
\F(u,\theta)=
\begin{pmatrix}
\dot{u}-f(.,u,\theta)\\
u(0) - u_0(\theta)\\
g(.,u,\theta)
  \end{pmatrix},
\end{align}
then the system can be written as the nonlinear operator equation
\begin{align*}
&\F(u,\theta)=\Y=(0,0,y). \quad\hphantom{hh}
\end{align*}
This fits into the Hilbert space framework by setting the function space for the state $u$ to
\begin{align}
\calU=H^1(0,T;V)
\end{align}
and $\calW, \calY$ are as in (\ref{spaceYW}). On the space $\calU$, we employ
\begin{align}
(u,v)_\calU=\int_0^T \paren{\dot{u}(t),\dot{v}(t)}_V dt + \paren{u(0),v(0)}_V \qquad
\end{align}
as an inner product, which induces a norm being equivalent to the standard norm \\$\sqrt{\int_0^T \|u(t)\|^2_V+\|\dot{u}(t)\|^2_V dt}$. This is the result of well-posedness in $\calU$ of the first order ODE $\dot{u}(t)=f(t), t \in (0,T), u(0)=u_0$ for $f \in L^2(0,T;V), u_0 \in V$.

The operator $\F$ is well-defined
by additionally imposing boundedness and local Lipschitz continuity on the functions inducing the Nemytskii operators. Differentiability of $\F$ also follows from the latter condition.

\begin{prop}	\label{All-at-onceDiff}
Let the Caratheodory mappings $f, g$ be:
\begin{enumerate}[label=(A\arabic*)]
\item
G\^ateaux differentiable with respect to their second and third arguments for almost all $t \in (0,T)$ with linear and continuous derivatives
\item \label{A2}
locally Lipschitz continuous in the sense
\begin{align*}
\forall M\geq 0, \exists L(M) &\geq 0, \forall^{a.e.} t \in (0,T):\\
& \|f(t,v_1,\theta_1)-f(t,v_2,\theta_2)\|_{V^*}
\leq L(M) (\|v_1-v_2\|_V +\|\theta_1-\theta_2\|_\calX),	\\
&\qquad\qquad\qquad\qquad \forall v_i \in V, \theta_i \in \calX: \|v_i\|_V, \|\theta_i\|_\calX \leq M, i=1,2
\end{align*}
and satisfy the boundedness condition, i.e., $f(.,0,0) \in \calW$.\\
The same is assumed to hold for $g$. Moreover, let $u_0$  also be G\^ateaux differentiable.
\end{enumerate}
Then $\F$ defined by (\ref{AAO-model}) is G\^ateaux differentiable on $\calU \times \calX$ and its derivative is given by 
\begin{align}	\label{AAO-adjoint}
&\F'(u,\theta)=
\begin{pmatrix}
(\frac{d}{dt}-f'_u(.,u,\theta)) & -f'_\theta(.,u,\theta)\\  
\delta_0 & -u'_0(\theta)\\
g'_u(.,u,\theta) & g'_\theta(.,u,\theta)
  \end{pmatrix}.
\end{align} 
\end{prop}

\proof
We show G\^ateaux differentiability of $f$ at an arbitrary point $(u,\theta) \in \calU \times \calX$. Without loss of generality, we consider the direction $(v,\xi)$ lying in a unit ball and $\epsilon \in (0,1]$
\begin{equation*}
\frac{1}{\epsilon}\| f(.,u+\epsilon v,\theta+\epsilon \xi) - f(.,u,\theta)
- \epsilon f'_u(.,u,\theta)v - \epsilon f'_\theta(.,u,\theta)\xi \|_\calW
= \paren{\int_0^T r_\epsilon(t)^2 dt}^\frac{1}{2},
\end{equation*}
where 
\begin{align} \label{r}
 r_\epsilon(t)&= \left\| \int_0^1 \bigg(
( f'_u(t,u(t)+\lambda\epsilon v(t),\theta +\lambda\epsilon \xi) - f'_u(t,u(t),\theta)) v(t) \right . \nonumber\\
&\qquad\qquad\qquad\left . + ( f'_\theta(t,u(t)+\lambda\epsilon v(t),\theta+\lambda\epsilon \xi) - f'_\theta(t,u(t),\theta)) \xi \bigg) d\lambda \right\|_{V^*}.
\end{align}
From local Lipschitz continuity and G\^ateaux differentiability of $f$, we deduce, by choosing $M=\max\{\sqrt{2T},\sqrt{2}\}(\|u\|_\calU +1) + \|\theta\|_\calX + 1$,
\begin{align*}
\|f'_u(t,u(t),\theta)w\|_{V^*}
&=\lim_{\epsilon \rightarrow 0}\left \|\frac{f(t,u(t)+\epsilon w,\theta)-f(t,u(t),\theta)}{\epsilon}\right\|_{V^*} \\
&\leq  L(M)\|w\|_V. 
\end{align*}
Continuity of the embedding $\calU\hookrightarrow C(0,T;V)$ is invoked above. Indeed, for any $t \in (0,T)$
\begin{flalign*}
&\quad \|u(t)\|_V 
\leq \int_0^T \|\dot{u}(s)\|_V ds+\|u(0)\|_V
\leq\max\{\sqrt{2T},\sqrt{2}\}\|u\|_\calU. \qquad
\end{flalign*}
As an immediate consequence
\begin{align*}
&\|f'_u(t,u(t)+\lambda\epsilon v(t),\theta+\lambda\epsilon\xi)\|_{V \rightarrow V^*} \leq L(M) \quad\\
&\|f'_\theta(t,u(t)+\lambda\epsilon v(t),\theta+\lambda\epsilon\xi)\|_{\calX \rightarrow V^*} \leq L(M) 
\end{align*}
for any $\lambda \leq 1, \epsilon\leq 1$. Substituting all the operator norm estimates into $r_\epsilon(t)$, we obtain
\begin{align}
r_\epsilon(t) \leq  2L(M)\paren{\|v(t)\|_V + \|\xi\|_\calX} := \bar{r}(t),
\end{align}
which is a square integrable function. Since $f$ is continuously differentiable, $r_\epsilon \rightarrow 0$ as $\epsilon\rightarrow 0$ for almost all every $t \in (0,T)$, and applying Lebesgue's Dominated Convergence Theorem yields G\^ateaux differentiability of $f$. Differentiability of $g$ is analogous.

In a similar manner, we can prove continuity of the derivative of $f, g$, also well-definedness of $f, g$ can be deduced from local Lipschitz continuity and the boundedness condition. \qed 

\begin{remark} \label{AAO-FrechetDiff}
If $f, g$ and $u_0$ are Fr\'echet differentiable then so is $\F$.
\end{remark}

\begin{remark}\text{ } 
\begin{itemize}
\item 
The local Lipschitz condition \ref{A2} is weaker than the one used in \cite{HaoDucThang}, as we only have the $V^*$-norm on the left hand side of the Lipschitz condition.
\item
Note that differentiablity of $\F$ can be interpreted on the stronger image space $C(0,T;V^*)$ since $r(.)$ is not only square integrable but also uniformly bounded with respect to time (provided that local Lipschitz continuity is fulfilled for all $t \in (0,T)$ instead of for almost all $t \in (0,T)$).
\item
Another idea could be a weakening of the Lipschitz continuity condition to
\begin{align*}
\forall M\geq 0, \exists L(M) &\geq 0, \gamma \in L^2(0,T), \forall^{a.e.} t \in (0,T):\\
& \|f(t,v_1,\theta_1)-f(t,v_2,\theta_2)\|_{V^*}
\leq L(M) \gamma(t)(\|v_1-v_2\|_V +\|\theta_1-\theta_2\|_\calX),	\\
&\qquad\qquad\qquad\qquad \forall v_i \in V, \theta_i \in \calX: \|v_i\|_V, \|\theta_i\|_\calX \leq M, i=1,2,
\end{align*}
then square integrability in time can be transferred from $\|v(.)\|_V$ to $\gamma(.)$.\\
\end{itemize}
\end{remark}

We now derive the Hilbert space adjoint for $\F'(u,\theta)$.
\begin{prop}	\label{All-at-once-Adjoint}
The Hilbert space adjoint of $\F'(u,\theta)$ is given by
\begin{flalign}
&\quad\quad\quad\quad \F'(u,\theta)^\star: \calW \times V \times \calY \rightarrow \calU \times \calX  \nonumber &\\
&\quad\quad\quad\quad \F'(u,\theta)^\star=
\begin{pmatrix}
(\frac{d}{dt}-f'_u(.,u,\theta))^\star & \delta_0^\star & g'_u(.,u,\theta)^\star\\
-f'_\theta(.,u,\theta)^\star & -u'_0(\theta)^\star & g'_\theta(.,u,\theta)^\star
  \end{pmatrix}, \label{All-at-onceop}&
\end{flalign}
where
\begin{flalign*}
&\quad\quad\quad\quad g'_u(.,u,\theta)^\star z = u^z &\\
&\quad\quad\quad\quad g'_\theta(.,u,\theta)^\star z = \int_0^T g'_\theta(t,u(t),\theta)^*z(t) dt &\\
&\quad\quad\quad\quad \delta_0^\star h= u^h &\\
&\quad\quad\quad\quad u'_0(\theta)^\star = u'_0(\theta)^*\tilde{I}^{-1} &\\	
&\quad\quad\quad\quad \big(\frac{d}{dt}-f'_u(.,u,\theta)\big)^\star w = u^w +\int_0^t \tilde{I}^*Iw(s)ds&\\
&\quad\quad\quad\quad f'_\theta(.,u,\theta)^\star w = \int_0^T f'_\theta(t,u(t),\theta)^*Iw(t) dt, &
\end{flalign*}
$u^z, u^w, u^h$ solve
\begin{flalign*}
&\quad\quad\quad\quad \begin{cases}
\ddot{u}^z(t)  = -\tilde{I}g'_u(t,u(t),\theta)^*z(t) \quad t \in(0,T)\\
\dot{u}^z(T)=0, \quad \dot{u}^z(0)- u^z(0)=0
\end{cases} &\\
&\quad\quad\quad\quad \begin{cases}
\ddot{u}^w(t) =  \tilde{I}f'_u(t,u(t),\theta)^*Iw(t) \quad t \in(0,T)\\
\dot{u}^w(T)=0, \quad \dot{u}^w(0)-u^w(0)=0
\end{cases} &\\
&\quad\quad\quad\quad \begin{cases}
\ddot{u}^h(t) = 0 \quad t \in(0,T)\\
\dot{u}^h(T)=0, \quad \dot{u}^h(0)-u^h(0)=-h,
\end{cases} &
\end{flalign*} 
and
\begin{flalign*}
&\quad\quad\quad\quad f'_u(t,u(t),\theta)^*: V^{**}=V \rightarrow V^*, \quad f'_\theta(t,u(t),\theta)^*: V^{**}=V \rightarrow \calX^*=\calX &\\
&\quad\quad\quad\quad g'_u(t,u(t),\theta)^*: Z^*=Z \rightarrow V^*, \quad\hphantom{f} g'_\theta(t,u(t),\theta)^*: Z^*=Z \rightarrow \calX^*=\calX &\\
&\quad\quad\quad\quad u'_0(\theta)^*:V^* \rightarrow \calX^*=\calX &
\end{flalign*}
are the respective Banach and Hilbert space adjoints.
\end{prop}

\proof

The adjoints of the derivatives with respect to $u$ can be established as follows.\\
For any $z \in \calY$ and $v \in \calU$
\begin{flalign*}
&\quad\paren{z,g'_u(.,u,\theta) v}_\calY &\\
&\quad\qquad= \int_0^T \paren{z(t),g'_u(t,u(t),\theta)v(t)}_Z dt &\\
&\quad\qquad = \int_0^T \dupair{g'_u(t,u(t),\theta)^*z(t),v(t)}_{V^*,V} dt&\\
&\quad\qquad = \int_0^T \paren{\tilde{I}g'_u(t,u(t),\theta)^*z(t),v(t)}_V dt&\\
&\quad\qquad = \int_0^T \paren{-\ddot{u}^z(t),v(t)}_V dt&\\
&\quad\qquad = \int_0^T \paren{\dot{u}^z(t),\dot{v}(t)}_V dt+\paren{u^z(0),v(0)}_V-\paren{\dot{u}^z(T),v(T)}_V+\paren{\dot{u}^z(0)-u^z(0),v(0)}_V&\\
&\quad\qquad = \paren{u^z,v}_\calU.&
\end{flalign*}
For any $w \in \calW$ and $v \in \calU$
\begin{flalign*}
&\quad\paren{w,\big(\frac{d}{dt}-f'_u(.,u,\theta)\big) v}_\calW &\\
&\quad\qquad= \int_0^T \paren{w(t),-f'_u(t,u(t),\theta) v(t)}_{V^*} dt +\int_0^T \paren{w(t),\frac{d}{dt} v(t)}_{V^*} dt&\\
&\quad\qquad= \int_0^T \dupair{-f'_u(t,u(t),\theta)^*Iw(t),v(t)}_{V^*,V} dt + \int_0^T \paren{Iw(t),\tilde{I}\frac{d}{dt}v(t)}_V dt&\\
&\quad\qquad= \int_0^T \paren{-\tilde{I}f'_u(t,u(t),\theta)^*Iw(t),v(t)}_V dt+\int_0^T \paren{\tilde{I}^*Iw(t),\frac{d}{dt}v(t)}_V dt&\\
&\quad\qquad= \int_0^T \paren{-\ddot{u}^w(t),v(t)}_V dt+\int_0^T \paren{\frac{d}{dt}\int_0^t \tilde{I}^*Iw(s)ds,\frac{d}{dt}v(t)}_V dt&\\
&\quad\qquad= \int_0^T \paren{\dot{u}^w(t),\dot{v}(t)}_V dt + \paren{u^w(0),v(0)}_V &\\
&\quad\qquad\qquad\qquad\qquad +\int_0^T \paren{\frac{d}{dt}\int_0^t \tilde{I}^*Iw(s)ds,\frac{d}{dt}v(t)}_V dt +\paren{\int_0^0 \tilde{I}^*Iw(s)ds,v(0)}_V&\\
&\quad\qquad= \paren{u^w+\int_0^t \tilde{I}^*Iw(s)ds,v}_\calU.&
\end{flalign*}
For any $h \in V$ and $v \in \calU$
\begin{flalign*}
&\quad\paren{h,\delta_0 v}_V 
= \int_0^T \paren{-\ddot{u}^h(t),v(t)}_V dt + \paren{h,v(0)}_V&\\
&\quad\qquad\qquad= \int_0^T \paren{\dot{u}^h(t),\dot{v}(t)}_V dt + \paren{u^h(0),v(0)}_V+ \paren{h+\dot{u}^h(0)-u^h(0),v(0)}_V&\\
&\quad\qquad\qquad= \paren{u^h,v}_\calU.&
\end{flalign*}
The three remaining adjoints are straightforward.
\qed

\begin{remark} \label{AAO-Kaczmarz}
For the Kaczmarz approach,
 the system is constructed as follows
\begin{align}
& \F_0(u,\theta)=
\begin{pmatrix}
\left(\dot{u}-f(.,u,\theta)\right)|_{(0,\tau_1)}\\
u(0) - u_0(\theta)\\
g(.,u,\theta)|_{(0,\tau_1)}
  \end{pmatrix}, 
& \F_j(u,\theta)=
\begin{pmatrix}
\left(\dot{u}-f(.,u,\theta)\right)|_{(\tau_j,\tau_{j+1})}\\
g(.,u,\theta)|_{(\tau_j,\tau_{j+1})}
  \end{pmatrix}\text{ } j=1\ldots n-1
\end{align} 
for $0=\tau_0<\tau_1<\ldots \tau_{n-1}<\tau_n=T$.
The function space setting for the Landweber-Kaczmarz method is
\begin{align*}
&\F_0: \calU \times \calX  \rightarrow \calW_0 \times V \times \calY_0, \qquad \F_j: \calU \times \calX  \rightarrow \calW_j \times \calY_j\quad j=1\ldots n-1
\end{align*} 
with
\begin{equation}
\calW_j=L^2(\tau_j,\tau_{j+1};V^*), \quad \calY_j=L^2(\tau_j,\tau_{j+1};Z)  \label{spaceWZj}
\end{equation}
thus
\begin{align}	
&\F'(u,\theta)_j^\star=
\begin{pmatrix}
\left(\left(\frac{d}{dt}-f'_u(.,u,\theta)\right)|_{(\tau_j,\tau_{j+1})}\right)^\star & \delta_0^\star & (g'_u(.,u,\theta)|_{(\tau_j,\tau_{j+1})})^\star\\
-(f'_\theta(.,u,\theta)|_{(\tau_j,\tau_{j+1})})^\star & -u'_0(\theta)^\star & (g'_\theta(.,u,\theta)|_{(\tau_j,\tau_{j+1})})^\star
  \end{pmatrix},
\end{align}
where the terms in the middle column $\delta_0^\star, -u'_0(\theta)^\star$ are present only in $\F'(u,\theta)_0^\star$.
\end{remark}

\medskip
\begin{remark} \label{AAO-adjoint-explicit}
In (\ref{All-at-onceop}), the adjoints of the derivatives with respect to $u$ can be represented explicitly
\begin{align*}
&\quad u^z(t)=  \int_0^T (t+1)\tilde{I}g'_u(s,u(s),\theta)^*z(s)ds - \int_0^t(t-s)\tilde{I}g'_u(s,u(s),\theta)^*z(s)ds \qquad\qquad\qquad\\
&\quad\tilde{u}^w(t):= u^w(t)+ \int_0^t \tilde{I}^*Iw(s)ds \\
&\quad\quad\quad \text{   }= -\int_0^T(t+1) \tilde{I}f'_u(s,u(s),\theta)^*Iw(s)ds + \int_0^t(t-s)\tilde{I}f'_u(s,u(s),\theta)^*Iw(s)ds \\
&\qquad\qquad\qquad\quad +\int_0^t \tilde{I}^*Iw(s)ds \\
&\quad u^h(t)= h 
\end{align*}
for $t \in (0,T).$

Incorporating the Kaczmarz scheme, we modify the adjoints accordingly
\begin{align*}
&\quad u^z(t)=  \int_{\tau_{j(k)}}^{\tau_{j(k+1)}} (t+1)\tilde{I}g'_u(s,u(s),\theta)^*z(s)ds - \int_{\min\{t,\tau_{j(k)}\}}^{\min\{t,\tau_{j(k+1)}\}}(t-s)\tilde{I}g'_u(s,u(s),\theta)^*z(s)ds \\
&\quad\tilde{u}^w(t)= -\int_{\tau_{j(k)}}^{\tau_{j(k+1)}}(t+1) \tilde{I}f'_u(s,u(s),\theta)^*Iw(s)ds\\
&\qquad\qquad\qquad\quad + \int_{\min\{t,\tau_{j(k)}\}}^{\min\{t,\tau_{j(k+1)}\}}(t-s)\tilde{I}f'_u(s,u(s),\theta)^*Iw(s)ds +\int_{\min\{t,\tau_{j(k)}\}}^{\min\{t,\tau_{j(k+1)}\}} \tilde{I}^*Iw(s)ds \\
&\quad u^h(t)= h
\end{align*}
for $t \in (0,T)$ with $j(k)=k-n\lfloor k/n\rfloor$, and also 
\begin{flalign*}
&\quad g'_\theta(.,u,\theta)^\star z = \int_{\tau_j}^{\tau_{j(k+1)}} g'_\theta(t,u(t),\theta)^*z(t) dt &\\
&\quad  u'_0(\theta)^\star = u'_0(\theta)^*I^{-1} &\\	
&\quad  f'_\theta(.,u,\theta)^\star w = \int_{\tau_j}^{\tau_{j(k+1)}} f'_\theta(t,u(t),\theta)^*Iw(t) dt.  &
\end{flalign*}
\end{remark}

\begin{remark} \label{AAO-discrete}
We now analyze the case of discrete measurement. Let $\{t_i\}_{i=1\ldots n}$ be the discrete observation time points, the system is now defined by
\begin{align}
&\F: \calU \times \calX  \rightarrow \calW \times V \times Z^n,
\qquad (u,\theta) \mapsto
\begin{pmatrix}
\left(\dot{u}-f(.,u,\theta)\right)\\
u(0) - u_0(\theta)\\
g_1(u,\theta)\\
\vdots\\
g_n(u,\theta)
  \end{pmatrix} \quad\hphantom{hh}
\end{align} 
with $g_i(u,\theta)=g(u,\theta)(t_i)=g(t_i,u(t_i),\theta), i=1\ldots n$, according to the definition of the Nemytskii operator. Differentiability of the Nemytskii operator $g_i$ could be inferred from differentiablity of the operator inducing it without the need of a local Lipschitz continuity condition.

The Hilbert space adjoint of $\F'(u,\theta)$ then takes the following form
\begin{align}
&\F'(u,\theta)^\star=
\begin{pmatrix}
(\frac{d}{dt}-f'_u(.,u,\theta))^\star & \delta_0^\star & g'_{iu}(u,\theta)^\star\\
-f'_\theta(.,u,\theta)^\star & -u'_0(\theta)^\star & g'_{i\theta}(u,\theta)^\star
  \end{pmatrix} \qquad\qquad\quad
\end{align}
in which the adjoint of $g'_u(u,\theta)$ mapping from the observation space $Z^n$ to $\calU$ is 
\begin{align*}
\paren{z,g'_u(u,\theta) v}_{Z^n}
&= \sum_{i=1}^n \paren{z_i,g'_{iu}(u,\theta)v(t_i)}_Z
= \sum_{i=1}^n \paren{\tilde{I}g'_{iu}(u,\theta)^*z_i,v(t_i)}_V \qquad\qquad\qquad\\\
&= \paren{u^z,v}_\calU,
\end{align*}
where
\begin{flalign*}
&\qquad\qquad\quad\begin{cases}
\ddot{u}^z(t) = 0 \quad t \in(0,T)\\
\dot{u}^z(0)-u^z(0)=0, \quad \dot{u}^z(t_i)=\tilde{I}g'_{iu}(u,\theta)^*z_i  \quad i=1\ldots n
\end{cases} &
\end{flalign*}
hence
\begin{align*}
u^z=\sum_{i=1}^n u^z_i
\qquad \text{with} \qquad
u^z_i=
\begin{cases}
  \tilde{I}g'_{iu}(u,\theta)^*z_i(t+1) \quad\text{ }  t \leq t_i\\
  \tilde{I}g'_{iu}(u,\theta)^*z_i(t_i+1)\quad  t > t_i.
\end{cases}
\end{align*}
If integrating into the  Kaczmarz scheme, on every subinterval of index $i$ we get $u^z=u^z_i$, so each equation in the system corresponds to one measurement.
\end{remark}

\color{blue}
\begin{remark}
Remarks \ref{AAO-adjoint-explicit} and \ref{AAO-discrete} show that the choice of the state space $\calU$, which is made in this paper, provides an explicit formula for the Hilbert space adjoint $\F'(u,\theta)^\star$. This enables us to speed up the computations in Landweber, Landweber-Kaczmarz and possibly in Newton-type methods.
\end{remark}
\color{black}
\section{Reduced setting} \label{secReduced}
In this section, we formulate the system into one operator mapping directly from the parameter space to the observation space. To this end, we introduce the parameter-to-state map
\[S:\calX\rightarrow\tilde{\calU}, \quad\text{ where }\quad u=S(\theta) \text{ solves } (\ref{model-1})-(\ref{model-2}),\]
then the forward operator for the Reduced setting can be expressed as
\begin{align} \label{Reducedop}
F:\calX \rightarrow \calY, \qquad \theta \mapsto g(.,S(\theta),\theta)
\end{align}
and the inverse problem of recovering $\theta$ from $y$ is
\begin{align*}
F(\theta)=y. \qquad
\end{align*}
Here, we use another state space $\tilde{\calU}$, which is different from $\calU$ in the All-at-once setting,
\begin{align}
\tilde{\calU}= W^{1,2,2}(0,T;V,V^*) \cap  L^\infty(0,T;V)
\end{align}
with $W^{1,2,2}(0,T;V,V^*)=\{u \in L^2(0,T;V): \dot{u} \in L^2(0,T;V^*)\}.$\\

To ensure existence of the parameter-to-state map, we assume that the operator $f$ meets the conditions in the following Assumption \ref{S} and the conditions \ref{R1} from Proposition \ref{reDiff}.

\begin{assumption} \label{S}
Let $\theta \in\calX$. Assume that
\begin{enumerate}[label=(S\arabic*)]
\item for almost all $t \in (0,T)$, the mapping $-f(t,.,\theta)$ is pseudomonotone \label{S1}
\item $ -f(.,.,\theta)$ is semi-coercive \label{S2}
\item $f$ satisfies a condition for uniqueness of the solution, e.g., 
\[ \forall u,v \in V, \forall^{a.e.}t \in (0,T): \langle f(t,u,\theta) - f(t,v,\theta),u - v\rangle_{V^*,V} \leq \rho^\theta(t)\|u -− v\|^2_H \]
for some $\rho^\theta \in L^1(0,T)$. \label{S3}
\end{enumerate}
\end{assumption}

\begin{prop} \label{reDiff}
Assume the model operator $f$ can be decomposed into the form
\begin{align*}
-f=f_1+f_2+f_3 \qquad\qquad
\end{align*}
with
\begin{flalign*}
&\qquad\qquad\hphantom{h} f_3:(0,T)\times \calX \rightarrow H  &\\
&\qquad\qquad\hphantom{h} f_2:(0,T)\times V \times \calX\rightarrow H &\\
&\qquad\qquad\hphantom{h} f_1:(0,T)\times V \times \calX \rightarrow V^*.
\quad &
\end{flalign*}
Additionally, let the following conditions be fulfilled
\begin{enumerate}[label=(R\arabic*)]
\item \label{R1}
the mappings $f_1, f_2, f_3$ satisfy
\begin{flalign*}	
&\qquad\hphantom{h} f_3 \in L^2(0,T;H) \text{ and continuous w.r.t. } \theta 	 &\\
&\qquad\hphantom{h}  \|f_2(t,v,\theta_\epsilon)\|_H \leq c_2^\theta(\gamma(t)+\|v\|_V) \text{ for some } \gamma \in L^2(0,T)  &\\
&\qquad\hphantom{h} f_1=\varphi'_v &\\
&\qquad\hphantom{h} \varphi:[0,T)\times V\times \calX \rightarrow \R \text{ convex w.r.t. $v$ and continuous}   &\\ 
&\qquad\hphantom{h} \varphi(t,v,\theta_\epsilon)\geq c_0^\theta\|v\|^2_V - c_1^\theta\|v\|^2_H     &\\
&\qquad\hphantom{h} \|\varphi'_t(t,v,\theta_\epsilon)\|_H \leq \tilde{c}_2^\theta(\tilde{\gamma}(t)+\|v\|_V^2) \text{ for some } \tilde{\gamma} \in L^1(0,T) &
\end{flalign*}
for all $\theta_\epsilon$ in some neighborhood of $\theta$ in $\calX$, with some $c^\theta_0>0,2c_1^\theta T < 1.$\\ If $\varphi:D(\varphi)\rightarrow \R^+, c^\theta_1$ does not need to be sufficiently small.
\item \label{R2}
$f$ is G\^ateaux differentiable with respect to its second and third arguments for almost all $t \in (0,T)$ with linear and continuous derivative.\\
\color{blue}
The derivative $-f'_u(.,u,\theta_\epsilon)$ moreover satisfies semi-coercivity in the sense
\[ \forall u, v \in V, \forall^{a.e.}t \in (0,T): \dupair{ -f'_u(t,u,\theta_\epsilon)v,v}_{V^*,V} \geq a_0^\theta|v|^2_V -a_1^\theta(t)|v|_V−-a_2^\theta(t)\|v\|_H^2 \]
for all $\theta_\epsilon$ in some neighborhood of $\theta$ in $\calX$, with some $a_0^\theta > 0, a_1^\theta \in L^2(0,T), a_2^\theta \in L^1(0,T)$ and some seminorm $|.|_V$ satisfying $\forall v \in V: \|v\|_V \leq c_{|.|}(|v|_V + \|v\|_H)$ for some $c_{|.|}>0.$
\color{black}
\item \label{R3}
$f$ is locally Lipschitz continuous in the sense
\begin{align*}
\forall M\geq 0, \exists L(M) &\geq 0, \forall^{a.e.} t \in (0,T):\\
& \|f(t,v_1,\theta_1)-f(t,v_2,\theta_2)\|_{V^*}
\leq L(M) (\|v_1-v_2\|_V +\|\theta_1-\theta_2\|_\calX),	\\
&\qquad\qquad\qquad\qquad \forall v_i \in V, \theta_i \in \calX: \|v_i\|_V, \|\theta_i\|_\calX \leq M, i=1,2.
\end{align*}
\end{enumerate}
Moreover, let also $u_0$ be G\^ateaux differentiable.
Then $F$, as defined by (\ref{Reducedop}), is G\^ateaux differentiable on $\calX$ and its derivative is given by 
\begin{align}
& F'(\theta) : \calX \rightarrow \calY \nonumber\\ 
& (F'(\theta)\xi)(t)=g'_u(t,S(\theta)(t), \theta)u^\xi(t)+g'_\theta(t,S(\theta)(t),\theta)\xi, \qquad\qquad \label{reF'}
\end{align}
where $u^\xi=S(\theta)'\xi$ solves the sensitivity equation
\begin{align}\label{S'eq}
&\begin{cases}
\dot{u}^\xi(t)=f'_u(t,S(\theta)(t),\theta)u^\xi(t)+f'_\theta(t,S(\theta)(t),\theta)\xi \quad t \in (0,T)\\ 
 u^\xi(0)=u'_0(\theta)\xi.
\end{cases}
\end{align}
\end{prop}

Before proving the result, we notice some facts.\\
Firstly, Assumption \ref{S} just guarantees that $S$ is a well-defined map from $\calX$ to $ W^{1,2,2}(0,T;V,V^*)$ \cite[Theorems 8.27, 8.31]{Roubicek}. To ensure it maps to $\tilde{\calU}= W^{1,2,2}(0,T;V,V^*)\cap L^\infty(0,T;V)$, condition \ref{R1} of Proposition \ref{reDiff} is additionally required. With condition \ref{R1}, we strengthen $S(\theta)$ to lie in $L^\infty(0,T;V)$. The proof is basically based on the regularity result in Roub\'i\v cek's book \cite{Roubicek} with extending the operator to be time-dependent (see the Appendix, Section \ref{Roubicek}).
 Secondly, due to the formulation, differentiability of the forward operator in the Reduced setting, in principle, is a question of differentiablity of $S$ (and $g$). \medskip

\textit{Proof of Proposition \ref{reDiff}.}
We show G\^ateaux differentiability of $S$. Fixing $\theta$ and without loss of generality, we consider $\xi$ lying in a unit ball and $\epsilon \in (0,1]$.\\
\hphantom{hhh} First, \ref{R1} allows us to apply the regularity result (\ref{regularity}) to obtain
\begin{align} \label{boundLinf}
\|S(\theta &+\epsilon\xi)\|_{L^\infty(0,T;V)}	\nonumber\\
&\leq N^\theta \bigg(2c^\theta_2\|\gamma\|_{L^2(0,T)} + 2\tilde{c}^\theta_2\|\tilde{\gamma}\|_{L^1(0,T)}^\frac{1}{2} + 2(c^\theta_2+\tilde{c}^\theta_2)\sqrt{\frac{c_1^\theta}{c_0^\theta}}\|u_0(\theta+\epsilon\xi)\|_H \nonumber\\
&\qquad\qquad\qquad\qquad\qquad + \|f_3(.,\theta+\epsilon\xi)\|_{L^2(0,T;H)} + \sqrt{|\varphi(0,u_0(\theta+\epsilon\xi),\theta+\epsilon\xi)|} \bigg) \nonumber\\
&\leq N^\theta \bigg( 2c^\theta_2\|\gamma\|_{L^2(0,T)}  + 2\tilde{c}^\theta_2\|\tilde{\gamma}\|_{L^1(0,T)}^\frac{1}{2} + 2C_{V\rightarrow H}(c^\theta_2+\tilde{c}^\theta_2)\sqrt{\frac{c_1^\theta}{c_0^\theta}}\|u_0(\theta)\|_V \nonumber\\
&\qquad\qquad\qquad\qquad\qquad  + \|f_3(.,\theta)\|_{L^2(0,T;H)} + \sqrt{|\varphi(0,u_0(\theta),\theta)|}+1 \bigg)	\\
&:= M^\theta \nonumber
\end{align}
for any $\epsilon \in [0,\bar{\epsilon}]$, where the constant $N^\theta$ depends only on $c^\theta_0,c^\theta_1,c^\theta_2, \tilde{c}^\theta_2, T$. Here we make use of continuity of the embedding $V\hookrightarrow H$ through the constant $C_{V\rightarrow H}$.

Let denote $v_\epsilon:=\frac{1}{\epsilon}\paren{S(\theta +\epsilon \xi) - S(\theta)}$. The function $v_\epsilon$ solves
\begin{align*}
\dot{v}_\epsilon(t)
&=\frac{1}{\epsilon}\paren{f(t,S(\theta +\epsilon\xi)(t),\theta +\epsilon\xi) - f(t,S(\theta)(t),\theta)}\\
&=\int_0^1 \paren{f'_u(t,S(\theta)(t) +\lambda\epsilon v_\epsilon(t),\theta +\lambda\epsilon\xi) v_\epsilon(t)
+f'_\theta(t,S(\theta)(t) +\lambda\epsilon v_\epsilon(t),\theta +\lambda\epsilon\xi)\xi}d\lambda\\
&=: A_\epsilon(t) v_\epsilon(t) + B_\epsilon(t)\xi\\
v_\epsilon(0)
&=\frac{1}{\epsilon}(u_0(\theta +\epsilon \xi) - u_0(\theta)),
\end{align*}
where by local Lipschitz continuity of $f$ and (\ref{boundLinf}) we have, for all most all $t \in (0,T),$
\begin{align*}
\|B_\epsilon(t)\|_{\calX \rightarrow V^*} &\leq L(\|S(\theta)\|_{L^\infty(0,T;V)}+\|S(\theta+\epsilon\xi)\|_{L^\infty(0,T;V)}+\|\theta\|_\calX+1) \qquad\qquad\quad\hphantom{h}\\
& \leq L(2M^\theta+\|\theta\|_\calX+1)\\
\|A_{\epsilon}(t)\|_{V\rightarrow V^*} &\leq  L( 2M^\theta+  \|\theta\|_\calX+1)
\end{align*}
for any $\epsilon \in (0,\bar{\epsilon}]$.\\
Together with {\ref{R2}}, we get the following estimate
\begin{align}
\|v_\epsilon\|_{W^{1,2,2}(0,T;V,V^*)} &\leq C^\theta \left( \|v_\epsilon(0)\|_H + \|B_\epsilon(.)\xi\|_{L^2(0,T;V^*)} \right) \qquad\qquad\quad\nonumber\\
&\leq C^\theta \left( C_{V \rightarrow H}\|u'_0(\theta)\|_{\calX \rightarrow V} + 1 + \sqrt{T} L(2M^\theta+\|\theta\|_\calX+1) \right) \quad\hphantom{.}
\end{align}
for any $\epsilon \in (0,\bar{\epsilon}].$

Let $\tilde{v}_\epsilon:=\frac{1}{\epsilon}\paren{S(\theta +\epsilon \xi) - S(\theta)-\epsilon u^\xi}$, then $\tilde{v}_\epsilon$ solves
\begin{align*}
\dot{\tilde{v}}_\epsilon(t)
&= f'_u(t,S(\theta)(t),\theta)\tilde{v}_\epsilon +\frac{1}{\epsilon} \bigg(-\epsilon f'_u(t,S(\theta)(t),\theta)v_\epsilon(t)-\epsilon f'_\theta(t,S(\theta)(t),\theta)\xi \qquad\qquad\\
&\qquad\qquad\qquad\qquad\qquad\quad  + f(t,S(\theta)(t) +\epsilon v_\epsilon(t),\theta +\epsilon\xi)-f(t,S(\theta)(t),\theta) \bigg)\\
&=: f'_u(t,S(\theta)(t),\theta)\tilde{v}_\epsilon + b_\epsilon(t)\\
\tilde{v}_\epsilon(0)&=\frac{1}{\epsilon}(u_0(\theta +\epsilon \xi)-u_0(\theta)-\epsilon u'_0(\theta)\xi ),
\end{align*}
as a result
\begin{equation} \label{Re-vtilde}
\|\tilde{v}_\epsilon\|_{W^{1,2,2}(0,T;V,V^*)} \leq  C^\theta \paren{\|\tilde{v}_\epsilon(0)\|_H+\|r_\epsilon\|_{L^2(0,T)}}.
\end{equation}
The fact that $v_\epsilon$ is bounded for any $\epsilon \in (0,\bar{\epsilon}]$ allows us to proceed analogously to the proof of Proposition \ref{All-at-onceDiff} to eventually get
\begin{equation*}
\|\tilde{v}_\epsilon\|_{W^{1,2,2}(0,T;V,V^*)} \rightarrow 0 \quad\text{as}\quad \epsilon \rightarrow 0,
\end{equation*}
which proves G\^ateaux differentiability of $S$.

Using the chain rule, we obtain the derivative of $F$ as in (\ref{reF'}) and the proof is complete.\qed

\begin{remark} \label{RED-FrechetDiff}
If $f, g$ and $u_0$ are Fr\'echet differentiable then so is $F$.
\end{remark}
Concluding differentiability of $F$, we now turn to derive the adjoint for $F'(\theta)$. 
\begin{prop}
The Hilbert space adjoint of $F'(\theta)$ is given by
\begin{align}
& F'(\theta)^\star  : \calY \rightarrow \calX \nonumber\\
& F'(\theta)^\star z= \int_0^T g'_\theta(t,S(\theta)(t), \theta)^* z(t) + f'_\theta(t,S(\theta)(t),\theta)^* p^z(t)dt + u'_0(\theta)^* p^z(0),
\end{align}
where $p^z$ solves 
\begin{align} \label{RED-pz}
&\begin{cases}
-\dot{p^z}(t)=f'_u(t,S(\theta)(t),\theta)^* p^z(t) + g'_u(t,S(\theta)(t),\theta)^* z(t) \quad t \in (0,T)\\
p^z(T)=0.
\end{cases}
\end{align}
\end{prop}
\proof
\cite[Proposition 2.7]{Kaltenbacher}. \qed

\begin{remark}
For the Kaczmarz approach and discrete measurement, we refer to \cite[Section 2.4, Remark 2.8]{Kaltenbacher} respectively.
\end{remark}

\color{blue}

\begin{remark}
With Propositions \ref{All-at-onceDiff} and \ref{reDiff} as well as Remarks \ref{AAO-FrechetDiff} and \ref{RED-FrechetDiff}, we obtain G\^ ateaux and Fr\' echet differentiability of the forward operators in both All-at-once and Reduced setting. Beyond the use of these derivatives in iterative methods (Landweber here, Gauss-Newton type methods in future work), knowledge of this differentiability yields more information on the topology of the function spaces. By utilizing Fr\' echet differentiability, other properties of the nonlinear operator can be understood.
\end{remark}
\section{Discussion} \label{secdiscuss}
\begin{remark}
In this paper, by introducing the new state spaces $\calU, \tilde{\calU}$ and imposing relevant structural conditions on the nonlinear forward operators, which mainly rely on local Lipschitz continuity, the restrictive growth conditions \cite[Conditions (2.2)-(2.5), (2.29)-(2.30)]{Kaltenbacher} in the All-at-once and Reduced formulations are eliminated.
\end{remark}
\subsection{Comparison between All-at-once and Reduced formulations}
Comparing the All-at-once and Reduced formulations of dynamic inverse problems, we observe the following
\begin{itemize}
\item
In the All-at-once setting, well-definedness of $\F$ directly follows from well-definedness of $f$ and $g$.
Meanwhile in the Reduced setting, it involves the need of well-definedness of the parameter-to-state map $S:\calX\rightarrow\tilde{\calU}$, which requires additional conditions on $f$ and $u_0$. Therefore, the All-at-once setting gives more flexibility in choosing the function spaces and can deal with more general classes of problems than the Reduced setting.
\item
In the general case, when $f: (0,T)\times V\times \calX \rightarrow W^*$ with the Hilbert space $W$ being possibly different from $V$, the local Lipschitz condition can be written as
\begin{align} \label{Lipschitzspace}
\forall M\geq 0, \exists L(M) &\geq 0, \forall^{a.e.} t \in (0,T):\\
& \|f(t,v_1,\theta_1)-f(t,v_2,\theta_2)\|_{W^*}
\leq L(M) (\|v_1-v_2\|_V +\|\theta_1-\theta_2\|_\calX),	\nonumber\\
&\qquad\qquad\qquad\qquad \forall v_i \in V, \theta_i \in \calX: \|v_i\|_V, \|\theta_i\|_\calX \leq M, i=1,2 \nonumber
\end{align}
and is applicable for both settings. With this condition, all the proofs proven in the All-at-once setting are unchanged, while in the Reduced setting, the proof for well-definedness of $S$ on the new function spaces might be altered.
\item 
The All-at-once version naturally carries over to the wave equation (or also fractional sub- or superdiffusion) context by just replacing the first time derivative by a
second (or fractional) time derivative. The Reduced version, however, requires additional conditions for well-definedness of the parameter-to-state map.
\end{itemize}
\subsection{On time-dependent parameter identification} \label{timedepend}
The parameter considered in the previous sections is time-independent. However, the theory developed in this paper, in principle, allows the case of time-dependent $\theta$ by introducing a time-dependent parameter space, for instance, $\calX=L^2(0,T;X), \calX=W^{1,2,2}(0,T;X,X^*)\times X_0, \calX=H^1(0,T;X)\times X_0$.

Relying on this function space, the local Lipschitz condition holding for All-at-once and Reduced setting reads as follows
\begin{align} \label{Lipschitztime}
\forall M\geq 0, \exists L(M) &\geq 0, \forall^{a.e.} t \in (0,T):\\
& \|f(t,v_1,\theta_1)-f(t,v_2,\theta_2)\|_{V^*}
\leq L(M) (\|v_1-v_2\|_V +\|\theta_1-\theta_2\|_X),	\nonumber\\
&\qquad\qquad\qquad\qquad \forall v_i \in V, \theta_i \in X: \|v_i\|_V, \|\theta_i\|_X \leq M, i=1,2. \nonumber
\end{align}
The proof for differentiability in both settings needs adapting to the new parameter function space. When working with concrete examples, if time point evaluation $\theta(t)$ is needed, the feasible choices are $\calX=H^1(0,T;X)\times X_0$ since $H^1(0,T)\hookrightarrow C(0,T)$ or $\calX=W^{1,2,2}(0,T;X,X^*)\times X_0$ since $W^{1,2,2}(0,T;X,X^*)\hookrightarrow C(0,T;H)$ with $X\hookrightarrow H\hookrightarrow X^*$. If the parameter plays the role of a source term, the time point evaluation will not be required.

In case $\calX=L^2(0,T;X)$ with Hilbert space $X$ and $u_0$ independent of $\theta$ (to avoid problems from existence of $\theta(0,.)$ for $\theta$ only in $L^2(0,T;X)$), the adjoint $\F'(u,\theta)^\star$ in the All-at-once setting is derived with the changes as follows
\begin{align}
f'_\theta(.,u,\theta)^\star w =  f'_\theta(.,u,\theta)^*Iw,
\qquad g'_\theta(.,u,\theta)^\star z =  g'_\theta(.,u,\theta)^*z  
\end{align}
with
\begin{align*}
f'_\theta(.,u,\theta)^*(t): V^{**}=V \rightarrow X^*=X, 
\qquad g'_\theta(.,u,\theta)^*(t): Z^*=Z \rightarrow X^*=X. 
\end{align*}
And for the the Reduced setting, one has
\begin{align}
F'(\theta)^\star z=  g'_\theta(.,S(\theta), \theta)^*z + f'_\theta(.,S(\theta),\theta)^* p^z,
\end{align}
where $p^z$ solves \eqref{RED-pz}.

The Kaczmarz approach relying on the idea of time segmenting does not directly carry over to this case and will be a subject for future research.
\color{black}
\section{The algorithm and its convergence} \label{secAlgorithm}
\subsection{Loping Landweber-Kaczmarz} \label{algorithm}
In this part, we write down the steps required in both settings of the Landweber-Kaczmarz method in case of continuous observation. Starting from an initial guess $(u_0, \theta_0)$, we run the iterations\medskip
\begin{paracol}{2}
\textit{All-at-once version}
\begin{enumerate}[label=\it S.\arabic*:]
\item Set argument to adjoint equations
\begin{align*}
& w_k(t) = \left(\dot{u}_k(t) - f(t,u_k(t),\theta_k)\right)|_{(\tau_j,\tau_{j+1})} \\ 
& h_k = u_k(0)-u_0(\theta_k)\\
& z_k(t) = \left(g(t,u_k(t),\theta_k)-y^\delta(t)\right)|_{(\tau_j,\tau_{j+1})} 
\end{align*}
so that 
\[(w_k,h_k,z_k)=\F_{j(k)}(u_k,\theta_k)-\Y_{j(k)}^\delta \qquad\]
\item Evaluate adjoint states
\begin{align*}
& \tilde{u}_k^w, f'_\theta(.,u_k,\theta_k)^\star w_k \quad\longleftarrow\quad w_k \qquad\qquad\qquad\\
&u_k^h, u'_k(0)^\star h_k \qquad\quad\hphantom{h} \longleftarrow\quad h_k
\end{align*}
\begin{align*}
&u_k^z, g'_\theta(.,u_k,\theta_k)^\star z_k \text{ }\quad\longleftarrow\quad z_k \qquad\qquad\qquad
\end{align*}
\item Update $(u,\theta)$ by 
\begin{align*}
& (u_{k+1},\theta_{k+1})= (u_k,\theta_k)\\
&\text{ } -\mu_k\F'_{j(k)}(u_k,\theta_k)^\star(\F_{j(k)}(u_k,\theta_k)-\Y_{j(k)}^\delta)
\end{align*}
\end{enumerate}
\switchcolumn
\textit{Reduced version}
\begin{enumerate}[label=\it S.\arabic*:]
\item Solve nonlinear state equation
\begin{align*}
&\begin{cases}
 \dot{u}_k(t) = f(t,u_k(t),\theta_k) \quad t\in (0,T) \qquad\qquad\qquad\text{ }\text{ } \\ 
 u_k(0)=u_0(\theta_k) 
\end{cases}
\end{align*}
\item Set argument to the adjoint equation
\begin{align*}
& z_k(t) = \left(g(t,u_k(t),\theta_k)-y^\delta(t)\right)|_{(\tau_j,\tau_{j+1})} \qquad
\end{align*}
so that \[z_k=F_{j(k)}(u_k,\theta_k)-y_{j(k)}^\delta \qquad\qquad\qquad\]
\item Evaluate the adjoint state 
\begin{align*}
& p_k^z \quad\longleftarrow\quad z_k \qquad\qquad\qquad\qquad\qquad\qquad
\end{align*}
\item Update $\theta$ by 
\begin{align*}
&\theta_{k+1}= \theta_k -\mu_k F'_{j(k)}(\theta_k)^\star(F_{j(k)}(\theta_k)-y_{j(k)}^\delta)
\end{align*} 
\end{enumerate}
\end{paracol}
with $j(k):=(k\mod n)$ for $n$ sub-intervals.

Now a stopping rule for the Landweber-Kaczmarz method should be specified. We illustrate the so-called loping strategy in the Reduced setting.
 We set
\begin{equation}
 \theta_{k+1}= \theta_k - w_k h_k,
\qquad\text{where}\qquad
h_k = \mu_kF'_{j(k)}(\theta_k)^\star(F_{j(k)}(\theta_k)-y_{j(k)}^\delta) 
\end{equation}
and
\begin{align} \label{w}
w_k=\begin{cases}
1 \qquad \text{if } \|F_{j(k)}(\theta_k)-y_{j(k)}^\delta\| \geq \tau\delta_{j(k)}\\
0 \qquad \text{otherwise}.
\end{cases}
\end{align}
Here $\tau>2$ is chosen subject to the tangential cone condition on $F$ (see below) and $\delta_{j(k)}$ is the noise level for the $j(k)$-th equation. In the Kaczmarz method, every collection of $n$ consecutive operators builds up a \enquote{cycle}. The iteration will process until it reaches the first index when the discrepancy principle holds over a full cycle. This means we stop the iteration at $k_*$, such that
\begin{align} \label{k}
w_{k_*-i}=0, \quad i=0,\ldots, n-1 \qquad\text{and}\qquad w_{k_*-n}=1.
\end{align}
The stopping rule \eqref{k} together with \eqref{w} corresponds to a discrepancy principle.

\begin{remark} \label{algorithm-observe}
We have some observations on the algorithm of the two settings.
\begin{itemize}
\item 
For each iteration, the All-at-once algorithm works only with linear models in all steps, while the Reduced algorithm requires one step solving a nonlinear equation to evaluate the parameter-to-state map. 
\item
Together with the fact that the adjoint states in the All-at-once setting can be analytically represented \textcolor{blue}{(see Remarks \ref{AAO-adjoint-explicit} and \ref{AAO-discrete})}, one step of the All-at-once algorithm is expected to run much faster than one of the Reduced algorithm.
\item
The residual in the All-at-once case comprises both the errors generated from $\theta$ and $u$ in the model and in the observations, while in the Reduced case, the exact $u=S(\theta)$ is supposed to be computed. Being inserted into the discrepancy principle, the stopping index $k_*$ in the Reduced algorithm is therefore possibly smaller than the one in the All-at-once case.
\item
For the Kaczmarz approach, the All-at-once setting restricts both model and observation operators to the subinterval $[\tau_j, \tau_{j+1}]$. The Reduced setting, however, applies the time restriction only for the observation operator, the model needs to be solved on the full time line to construct the parameter-to-state map.
\item
One can also incorporate the Kaczmarz strategy into the discrete observation case for both settings.  
\end{itemize}
\end{remark}

\subsection{Convergence analysis}
We are now in a position to provide convergence results under certain conditions. These conditions are derived in the context of iterative regularization methods for nonlinear ill-posed problems \cite{BarbaraNeubauerScherzer}.

\begin{assumption}\text{}\label{convergenceAss}
\begin{itemize}
\item Tangential cone condition in the All-at-once version
\begin{align} \label{AAO-TCC}
&\|f(.,\tilde{u},\tilde{\theta}) - f(.,u,\theta) - f'_u(.,u,\theta)(\tilde{u}-u) - f'_\theta(.,u,\theta)(\tilde{\theta}-\theta) \|_{\calW} \qquad\qquad\qquad\qquad \nonumber\\
&\qquad+ \|u_0(\tilde{\theta}) - u_0(\theta) - u'_0(\theta)(\tilde{\theta}-\theta) \|_V \\
&\qquad + \|g(.,\tilde{u},\tilde{\theta}) - g(.,u,\theta) - g'_u(.,u,\theta)(\tilde{u}-u) - g'_\theta(.,u,\theta)(\tilde{\theta}-\theta) \|_{\calY} \nonumber\\
& \leq c_{tc} \left( \| \dot{\tilde{u}}-\dot{u} - f(.,\tilde{u},\tilde{\theta}) + f(.,u,\theta) \|_{\calW} + \|u_0(\tilde{\theta}) - u_0(\theta)\|_V+ \|g(.,\tilde{u},\tilde{\theta}) - g(.,u,\theta) \|_{\calY} \right) \nonumber
\end{align}
and in the Reduced version
\begin{align} \label{Re-TCC}
&\|g(.,S(\tilde{\theta}),\tilde{\theta}) - g(.,S(\theta),\theta) - g'_u(.,S(\theta),\theta)v - g'_\theta(.,S(\theta),\theta)(\tilde{\theta}-\theta) \|_{\calY} \quad\qquad\qquad\qquad\qquad\qquad\nonumber\\
& \leq \tilde{c}_{tc}\|g(.,S(\tilde{\theta}),\tilde{\theta}) - g(.,S(\theta),\theta) \|_{\calY}
\end{align}
for some $c_{tc}, \tilde{c}_{tc}<\frac{1}{2}$, where $v$ solves
\begin{flalign} \label{Re-TTCv}
\quad &\begin{cases}
\dot{v}(t)=f'_u(t,S(\theta)(t),\theta)v(t)+f'_\theta(t,S(\theta)(t),\theta)(\tilde{\theta}-\theta) \quad t \in (0,T)\\
v(0)=u'_0(\theta)(\tilde{\theta}-\theta)
\end{cases}&
\end{flalign}

\item
The constant in the discrepancy principle is sufficiently large, i.e., $\tau > 2\dfrac{1+c_{tc}}{1-2c_{tc}}$

\item
The stepsize parameter satisfies $\mu_k \in \Big(0,\dfrac{1}{\| \F'(u_k,\theta_k)\|^2}\Big].$
\end{itemize}
\end{assumption}

Since our methods are considered in Hilbert spaces, we can employ existing convergence results collected from the book \textcolor{blue}{\cite[Theorem 3.26]{BarbaraNeubauerScherzer}}.

\begin{coro}
Let the assumptions of Proposition \ref{All-at-onceDiff} and Assumption \ref{convergenceAss} hold. Moreover, let the stopping index $k_*=k_*(\delta,y^\delta)$ be chosen as in (\ref{k}). Then the Landweber-Kaczmarz iteration in the All-at-once setting converges to a solution of $\F(u^\dagger,\theta^\dagger)=\Y$, provided that the starting point $(u_0,\theta_0)$ is sufficiently close to $(u^\dagger,\theta^\dagger)$.
\end{coro}
\begin{coro}
Let the assumptions of Proposition \ref{reDiff} and  Assumption \ref{convergenceAss} hold. Moreover, let the stopping index $k_*=k_*(\delta,y^\delta)$ be chosen as in (\ref{k}). Then the Landweber-Kaczmarz iteration in the Reduced setting converges to a solution of \hphantom{hh} $F(\theta^\dagger)=y$, provided that the starting point $\theta_0$ is sufficiently close to $\theta^\dagger$.
\end{coro}

\section{An example} \label{secEx}
In this section, the first part is dedicated to examining the conditions proposed in the abstract theory for a class of problems. This work will expose the maximum nonlinearity allowed in our setting, which indicates the improvement comparing to the current result \cite{Kaltenbacher}. The section is subsequently continued by some numerical experiments running on both continuous and discrete observations.

\subsection{Model for a class of problem}
Let us consider the semilinear diffusion system
\begin{alignat}{3}
&\dot{u}=\Delta u -\Phi(u) + \theta \qquad &&(t,x) \in (0,T)\times \Omega& \label{ex1}\\
& u_{|\partial\Omega}=0 && t \in (0,T)&	\label{ex2}\\
& u(0)=u_0 && x \in \Omega &	\label{ex3}\\
& y=Cu && (t,x) \in (0,T)\times \Omega,	\label{ex4}&
\end{alignat}
where $\Omega \subset \R^d$ is a bounded Lipschitz domain.\\
We investigate this  problem in the function spaces
\begin{flalign*} 
\calX=L^2(\Omega),  \qquad  V=H_0^1(\Omega), \qquad H=L^2(\Omega), \qquad Z=L^2(\Omega) 
\end{flalign*}
with linear observation (i.e., $C$ is a linear operator).
 For the Reduced setting, we assume that the nonlinear operator $\Phi$ is monotone and $\Phi(0)=0$.
By this way, one typical example could be given, for instance, $\Phi(u)=|u|^{\gamma-1}u, \gamma\geq 1.$

We now verify the imposed conditions in the All-at-once and Reduced versions. To begin, we decompose the model operator into
\[-f = -\Delta u +\Phi(u) - \theta := f_1 + f_3. \qquad\]
It is obvious that $f_3=-\theta \in L^2(0,T;H)$, and $g$ is Lipschitz continuous with the Lipschitz constant $L(M)=\|C\|_{V\rightarrow Z}$ and satisfies the boundedness condition $g(t,0,0)=0 \in Z$.

The next part
focuses on analyzing the properties of the model operator $f$. 

\begin{enumerate}[label=(R\arabic*)]
\item \textit{Regularity conditions for $f_1$}\medbreak
The regularity condition \ref{R1} holds by the following argument
\begin{align*}
& f_1(v)=-\Delta v +\Phi(v) = \phi'(v) \\
&\phi(v) = \int_\Omega \frac{1}{2}|\nabla v|^2 + \frac{1}{\gamma+1}\Phi(v)v dx \\
&\qquad \geq \frac{1}{2}\|v\|^2_V \text{ and continuous.}\\
&\phi''(v)[w,w] = \int_\Omega |\nabla w|^2 + \Phi'(v)ww dx \geq 0, \quad \forall v, w \in V \qquad\qquad\qquad\qquad
\end{align*}
concludes convexity of $\phi$. Here, we invoke monotonicity and differentiability of $\Phi$.\\
Analogously, $-f'_u(.,u,\theta_\epsilon)=-\Delta +\Phi'(u)$ is semi-coercive in the sence of \ref{R2}.
\end{enumerate}

\begin{enumerate}[label=(R\arabic*)]
\setcounter{enumi}{2}
\item \textit{Local Lipschitz continuity of $f$}\medbreak
First, we observe that
\begin{align*}
& \|f(t,v_1,\theta_1)-f(t,v_2,\theta_2)\|_{V^*} = \|\Phi(v_2)-\Phi(v_1) - \Delta(v_2-v_1)-(\theta_2-\theta_1) \|_{H^{-1}} \quad\\
&\qquad= \sup_{ \|\nabla w\|_{L^2(\Omega)} \leq 1} \int_\Omega (\Phi(v_2)-\Phi(v_1))w + \nabla(v_2-v_1)\nabla w - (\theta_2-\theta_1)w dx \\
&\qquad \leq \sup_{ \|\nabla w\|_{L^2(\Omega)} \leq 1} \int_\Omega (\Phi(v_2)-\Phi(v_1))w dx + \|v_1-v_2\|_V + \|\theta_1-\theta_2\|_{\calX}.c_{PF}, 
\end{align*}
where $c_{PF}$ is the constant in the Poincar\'e-Friedrichs inequality:
$\|v\|_{L^2(\Omega)} \leq c_{PF}\|\nabla v\|_{L^2(\Omega)}, \forall v \in H_0^1(\Omega).$
Developing the first term on the right hand side by applying H\"older's inequality, we have
\begin{flalign*}
&\sup_{ \|\nabla w\|_{L^2(\Omega)}  \leq 1} \int_\Omega (\Phi(v_2)-\Phi(v_1))w dx = \sup_{ \|\nabla w\|_{L^2(\Omega)} \leq 1} \left( \int_\Omega w^{\bar{p}}dx \right)^\frac{1}{\bar{p}} \left( \int_\Omega (\Phi(v_2)-\Phi(v_1))^\frac{\bar{p}}{\bar{p}-1}dx \right)^\frac{\bar{p}-1}{\bar{p}} &\\
&\quad \leq c_{H^1\rightarrow L^{\bar{p}}}\sqrt{1+c^2_{PF}}\left[ \int_\Omega (v_1-v_2)^{\bar{p}}dx \right]^\frac{1}{\bar{p}}. \left[ \int_\Omega \left( \int_0^1\Phi'(v_1+\lambda(v_2-v_1))d\lambda \right)^\frac{\bar{p}}{\bar{p}-2} dx \right]^\frac{\bar{p}-2}{\bar{p}} &\\
&\quad \leq \gamma c_\gamma c^2_{H^1\rightarrow L^{\bar{p}}}(1+c^2_{PF})\|v_1-v_2\|_V \left[ \int_0^1 \int_\Omega \left( \lambda v_2+(1-\lambda)v_1 \right)^{(\gamma-1)\frac{\bar{p}}{\bar{p}-2}}dx d\lambda \right]^\frac{\bar{p}-2}{\bar{p}} &\\
&\quad \leq 2^{\gamma-1}\gamma c_\gamma c^2_{H^1\rightarrow L^{\bar{p}}}(1+c^2_{PF})\|v_1-v_2\|_V \left( \|v_1^{\gamma-1}\|_{L^\frac{\bar{p}}{\bar{p}-2}(\Omega)} +\|v_2^{\gamma-1}\|_{L^\frac{\bar{p}}{\bar{p}-2}(\Omega)} \right), &
\end{flalign*}
provided additionally that $|\Phi'(v)|\leq c_\gamma|v|^{\gamma-1}.$
Altogether, we arrive at
\begin{align*}
& \|f(t,v_1,\theta_1)-f(t,v_2,\theta_2)\|_{V^*} \leq L(\|v_1\|^{\gamma-1}_V+\|v_2\|^{\gamma-1}_V)(\|v_1-v_2\|_V + \|\theta_1-\theta_2\|_\calX)
\end{align*}
if
\begin{align*}
(\gamma-1)\frac{\bar{p}}{\bar{p}-2} \leq \bar{p} \quad\Leftrightarrow\quad \gamma \leq \bar{p}-1,
\end{align*}
where $c_{H^1\rightarrow L^{\bar{p}}}$ is the constant and $\bar{p}$ is the maximum power such that the embedding $H^1(\Omega) \hookrightarrow L^{\bar{p}}(\Omega)$ is continuous
\begin{align}	\label{cond-power}
\bar{p}
\begin{cases}
=\infty				\quad &\text{if } d=1\\
<\infty 			\quad &\text{if } d=2\\
=\frac{2d}{d-2}		\quad &\text{if } d\geq 3 
\end{cases}
\quad \Leftrightarrow\quad
\gamma
\begin{cases}
=\infty				\quad &\text{if } d=1\\
<\infty 			\quad &\text{if } d=2\\
\leq 5					\quad &\text{if } d=3\\
\leq \frac{d+2}{d-2}	\quad &\text{if } d\geq 4. 
\end{cases} \qquad
\end{align}
\end{enumerate}
\color{blue}
\begin{remark}
Condition \eqref{cond-power} on $\gamma$ reveals that in two dimensional space, our proposed setting works out well with all finite powers of the nonlinearity and in one dimensional space, the accepted power is unconstrained. Noticeably, in the practical case, i.e., three dimensions, the largest power we attain is up to $5$ which enhances the limit in the nonlinearity of the current work and enables us to include important applications that had been ruled out in \cite{Kaltenbacher} due to the growth conditions there.\\
\end{remark}
\color{black}
\quad The remaining task is to examine well-definedness of the parameter-to-state map $S:\calX \rightarrow W^{1,2,2}(0,T;V,V^*)$, which is straightforward and the tangential cone condition for the forward operators, which is presented in Section \ref{TCC} in Appendix. For verifying the tangential cone condition, we additionally impose the growth condition $|\Phi''(v)|\leq c_\gamma|v|^{\gamma-2}$.

\begin{remark} \label{uniqueness}
At this point, we can briefly discuss the uniqueness question.
Since $\theta$ is time-independent, if having observation at any single time instance $t \in (0,T)$ and \textcolor{blue}{observation on all of $\Omega$}, we are able to determine $\theta$ uniquely by
\begin{align}\label{directinversion}
\theta=\left(\dot{u}-\Delta u +\Phi(u)\right)(t).
\end{align} 
The point evaluation at $t$ only works for sufficiently smooth observation, e.g., $u \in C^1(I;L^2(\Omega)) \cap C(I;H^2(\Omega))$ for some neighborhood $I$ of $t$, so that $\theta \in L^2(\Omega)$. \textcolor{blue}{In principle, this also induces a reconstruction scheme, namely, after filtering the given noisy data, applying formula \eqref{directinversion}. However, in contrast to the scheme we propose here, this would not apply to the practically relevant case of partial (e.g. boundary) and time discrete observations.}\\
Alternatively, we can approach this issue by using Theorem 2.4 \cite{BarbaraNeubauerScherzer}, which concludes uniqueness of the solution from the tangential cone condition as in Section \ref{secAlgorithm} and the null space condition
\[ \calN(F'(\theta^\dagger)) \subset \calN(F'(\theta)) \quad \text{for all } \theta \in \calB_\rho(\theta^\dagger), \qquad\]  
where $\calB_\rho(\theta^\dagger)$ is a closed ball of radius $\rho$ centering at $\theta^\dagger$. Since the tangential cone condition has been verified, it remains to examine the null space of $F'(\theta^\dagger)$. This could be done by linearizing the equation then utilizing some results from Isakov's book \cite{Isakov}.\\
The null space of $F'(\theta^\dagger)$ consists all $\theta$ such that the solution $v$ of
\begin{flalign*}
&\qquad\qquad\qquad\dot{v}-\Delta v +\Phi'(u^\dagger)v = \theta \qquad\hphantom{l} (t,x) \in (0,T)\times \Omega & \\
&\qquad\qquad\qquad v_{|\partial\Omega}=0 \qquad\qquad\qquad\qquad t \in (0,T)&	\\
&\qquad\qquad\qquad v(0)=0 \qquad\qquad\qquad\quad\hphantom{hi} x \in \Omega &
\end{flalign*}
leads to vanishing observations.\\
Theorem 9.1.2 \cite{Isakov} states that the solution to this inverse source problem is unique if a final observation is available and $\theta$ has compact support. As a result, the observation operator just needs to convey the information of $u$, e.g., at the final time $g(u)=u(T)$.\\
In case of discrete measurement and $\theta$ does not have compact support, uniqueness is still attainable by using Theorem 9.2.6 \cite{Isakov}. This theorem supports the case when one can measure one intermediate time data, i.e., $u(t_i,x) \text{ for some } t_i\in(0,T)$. In this situation, also observation of Neumann data on an arbitrary part of the boundary is demanded. 
\end{remark}

\subsection{Numerical experiment}
In the following numerical experiment, we select the nonlinear term $\Phi(u)= u^3$ motivated by the superconductivity example. We assume to observe  $u$ fully in time and space, i.e., $Cu=u \text{ on } (0,T)\times\Omega $ and that at initial time $u(0)=u_0=0$. The method in use is loping Landweber-Kaczmarz.

The parameters for implementation are as follows: the time line (0,0.1) (101 time steps) is segmented into 5 time subintervals, the space domain is $\Omega=(0,1)$ (101 grid points) and the system is perturbed by $5\%$ data noise.

\begin{figure}[!t] 
\setlength{\belowcaptionskip}{-5pt}
\centering
\includegraphics[width=4.5cm]{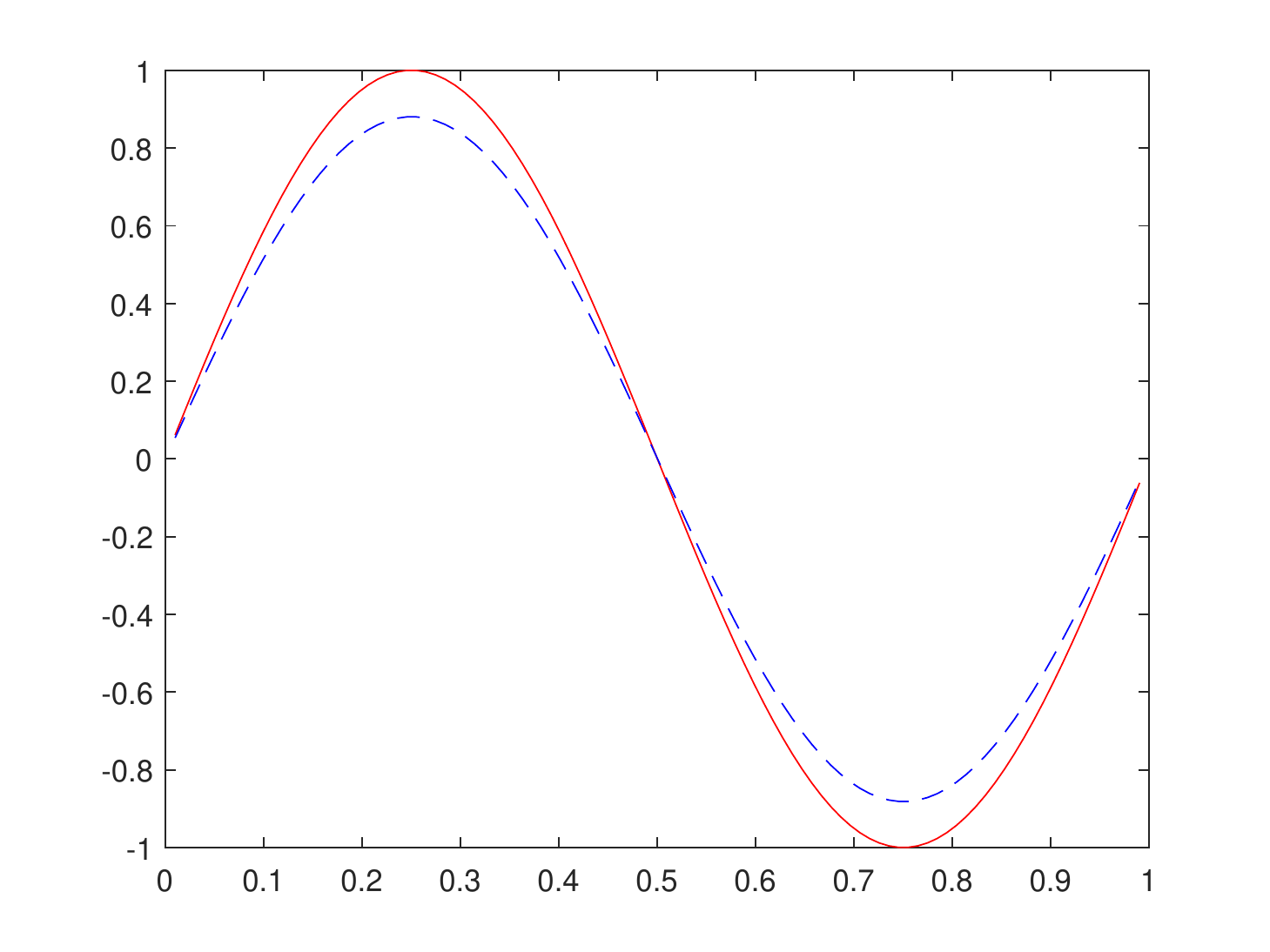}
\includegraphics[width=4.5cm]{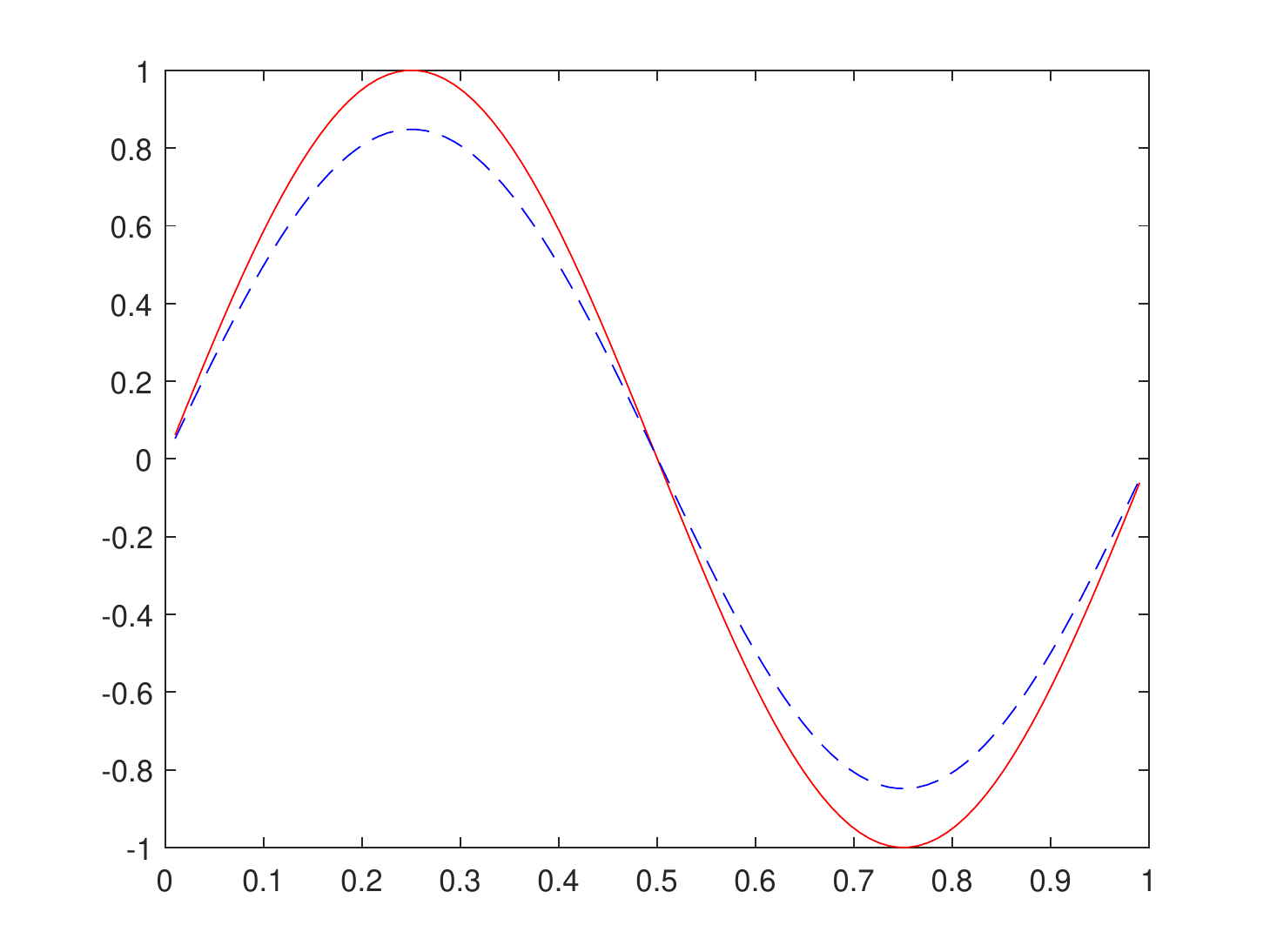}\\\vspace{-10pt}
\includegraphics[width=4.5cm]{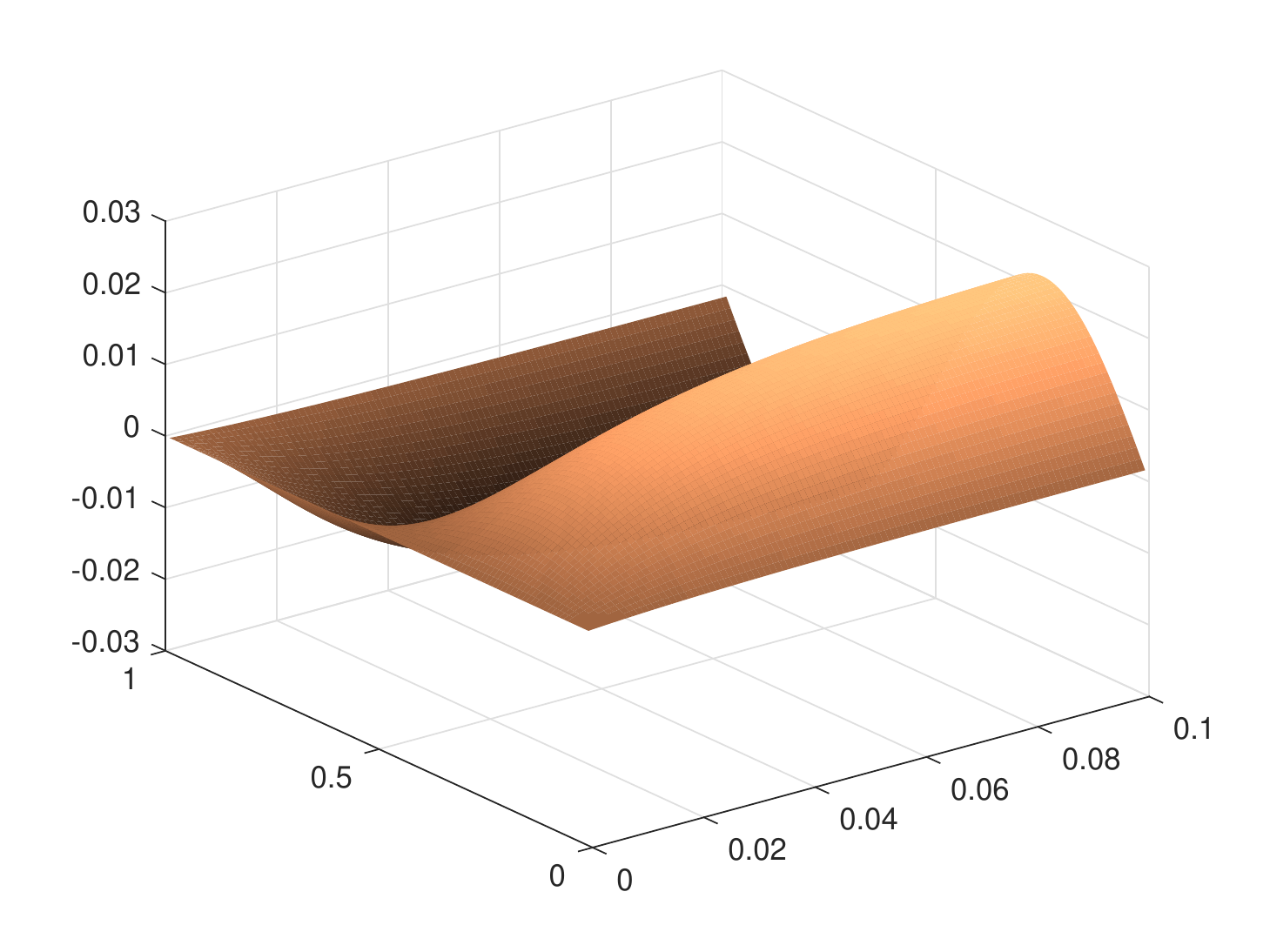}
\includegraphics[width=4.5cm]{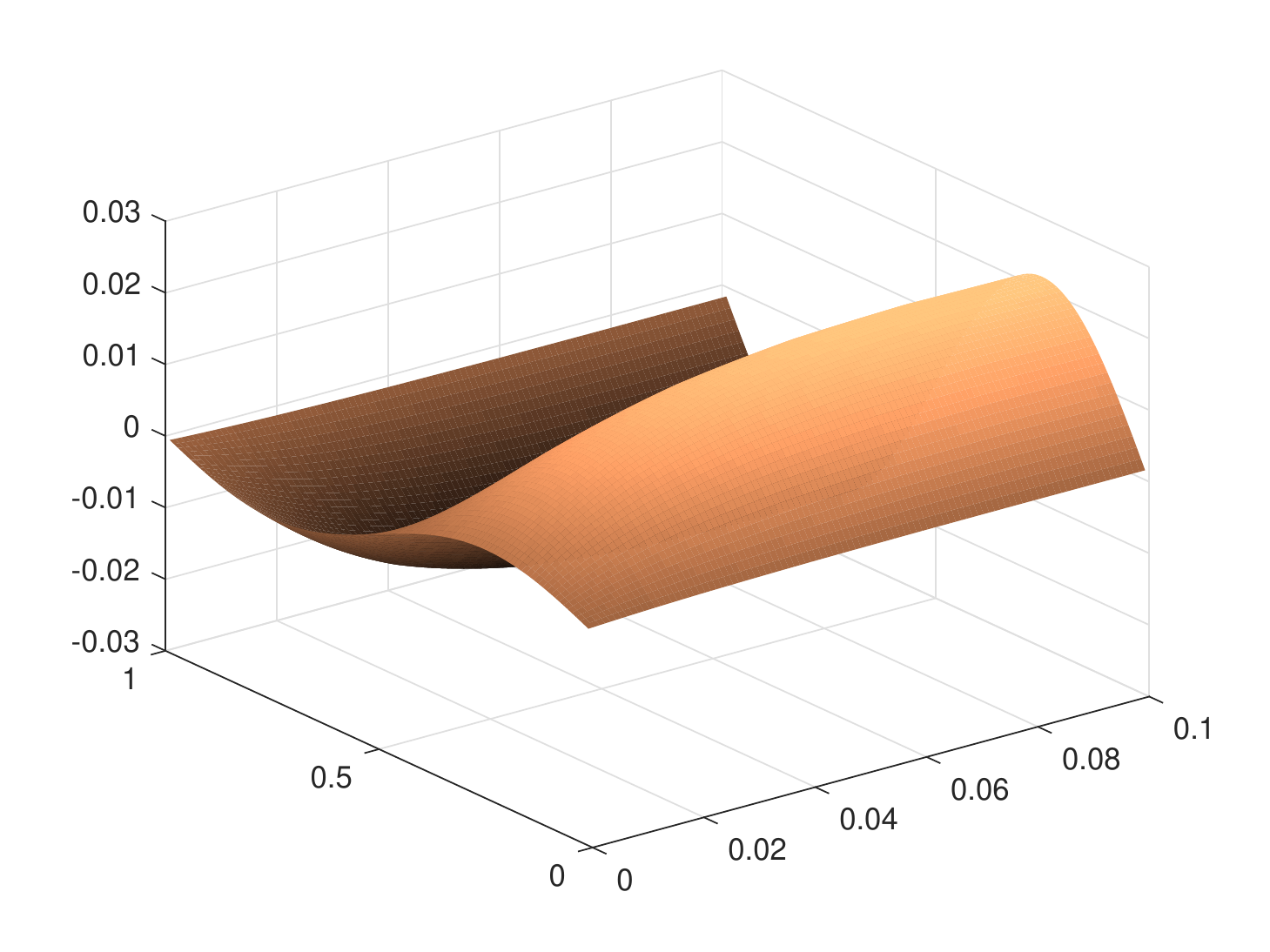}\\\vspace{-10pt}
\includegraphics[width=4.5cm]{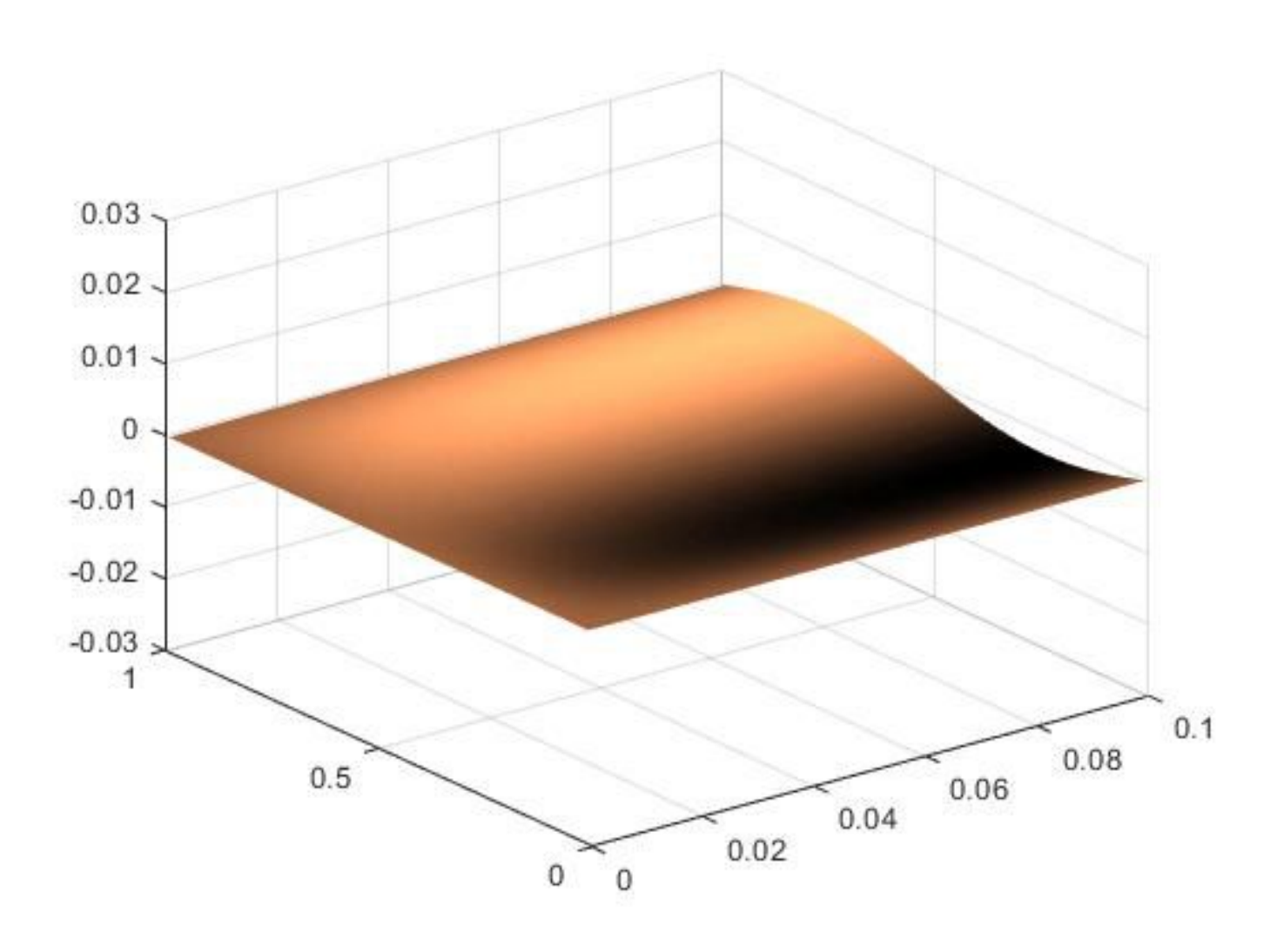}
\includegraphics[width=4.5cm]{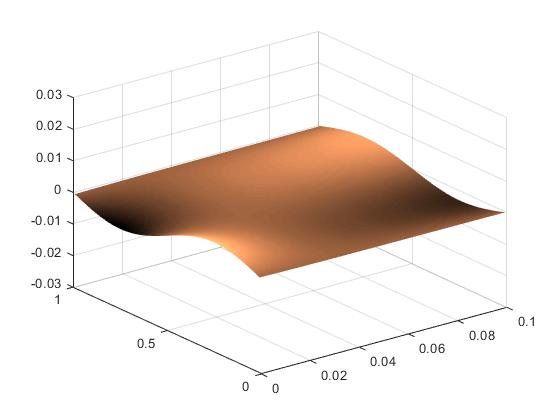}
\caption{Top: exact parameter $\theta^\dagger$ (solid) and reconstruction (dashed). Middle: reconstructed state $u$. Bottom: difference $u-u^\dagger$. Left: Reduced, right: AAO setting.}
\label{fig-1}
\end{figure}

Figure \ref{fig-1} displays the results of reconstructed parameter and state comparing to the exact ones. Apparently, two settings yield very similar results, except at $t=0$ where the Reduced setting approximates the exact initial state better. This is explained by the fact that the model equations (\ref{model-1})-(\ref{model-2}) in the Reduced setting is fully preserved to construct the parameter-to-state map, while in the All-at-once (AAO) setting, $u_0$ only appears in the forward operator with the index zero. However, the reconstructed parameters in both settings are definitely comparable.

\begin{figure}[!tb] 
\centering
\includegraphics[width=4.3cm]{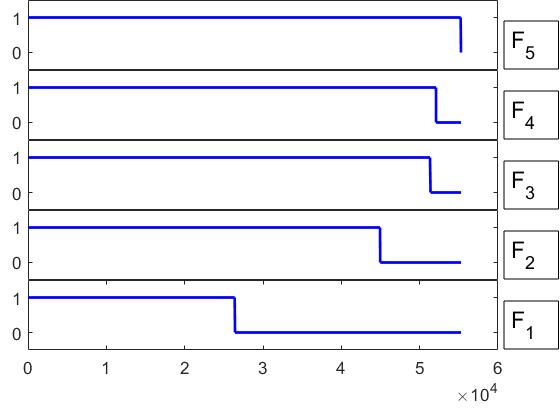}\quad
\includegraphics[width=4.3cm]{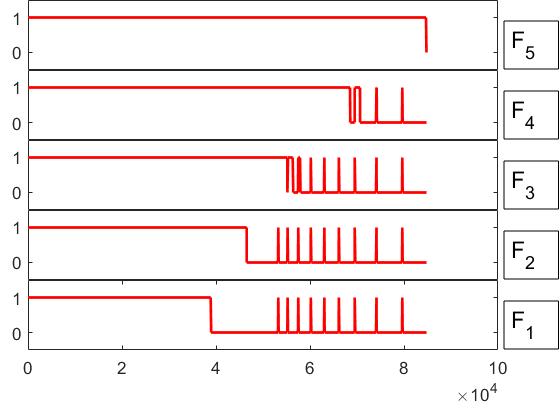}\quad
\includegraphics[width=4.3cm]{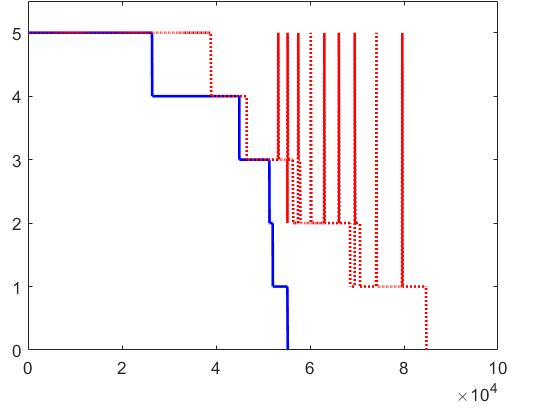}\quad
\caption{On/off iterations on each sub-interval per cycle for Reduced (left), AAO (middle) methods. Right: Number of inner steps per cycle for Reduced (solid blue), AAO (dashed red) settings. }
\label{fig-2}
\end{figure}

In Figure \ref{fig-2}, the left and the middle figures show the scalar $w_k$ in (\ref{w}) for each iteration (horizontal axis). Five operators represent five time subintervals correspondingly. The right figure sums all $w_k$ over five subintervals and plots both settings together. The All-at-once setting stops the iterations after a factor of $1.5$ times those in the Reduced setting. This means the All-at-once setting requires much more iterations than the Reduced one to obtain the accepted error level. Nevertheless, it runs much faster than the Reduced case, in particular the cpu times are: 2989s (Reduced method), 1175s (AAO method). The reasons for this effect have been discussed in Remark \ref{algorithm-observe}.

To demonstrate the discrete observation case, we ran tests at several numbers of discrete observation time points. The parameters for implementation stay the same as in the continuous observation case. Here we also use the loping Landweber-Kaczmarz method.

Table \ref{tab-1} compares the Reduced with the All-at-once setting at different numbers of observation time points $np$, where each of the time points corresponds to a sub equation, i.e., $np=n$ in the Kaczmarz method (\ref{KaczmarzIntro}). Those numbers vary largely from $3$ to the maximum discretization time step $101$. Despite the largely varying $np$, the errors in both settings are quite stable. It gets plausible when looking at the \enquote{\#updates} columns reporting the \enquote{real} iterations, which are the iterations making the update, i.e., iterations with $w_k=1$. This reveals that significantly increasing the number of discrete observation time points does not bring any improvement. We reason this phenomenon by the above argument (cf. Remark \ref{uniqueness}) of uniqueness of the solution according to which, one additional intermediate time data is sufficient to uniquely recover $\theta$.

\begin{table}[!t] 
\caption{Numerical experiment with different numbers of discrete observation points (np) at 5\% noise. Observation points are uniformly distributed on $(0,T].$}
\centering
\begin{tabular}{  c | R{1.3cm} R{1.3cm} c c | R{1.3cm} R{1.3cm} c c  }
\hline
np& &  \text{ }Reduced & & & & \hspace{3pt}AAO &\\
\hline
 & \#loops & \#updates & \hphantom{hh}cpu(s) & $\frac{\|\theta^\delta_{k_*}-\theta^\dagger\|}{\|\theta^\dagger\|}$ &  \#loops & \#updates & \hphantom{hh}cpu(s) & $\frac{\|\theta^\delta_{k_*}-\theta^\dagger\|}{\|\theta^\dagger\|}$ \\ 
\hline
\rule{0pt}{2ex}3 
&  6191&  3854 & 38 &   0.113 
&  9722&  6215 & 17 &   0.158 \\ 
\rule{0pt}{2ex}11 
&  6214&  4506 & 39 &   0.113 
&  9921&  7142 & 18 &   0.140 \\ 
\rule{0pt}{2ex}21 
&  6215&  4539 & 39 &   0.113 
&  10793&  7253 & 19 &   0.137 \\ 
\rule{0pt}{2ex}51 
&  6278&  4556 & 42 &   0.112 
&  11831&  7303 & 20 &   0.135 \\ 
\rule{0pt}{2ex}101 
&  9896 &  4626 & 62 &   0.108 
&  16865&  7389 & 28 &   0.132 \\ 
\hline
\end{tabular}
\label{tab-1}
\end{table}

\begin{figure}[!htb] 
\centering
\includegraphics[width=3.7cm]{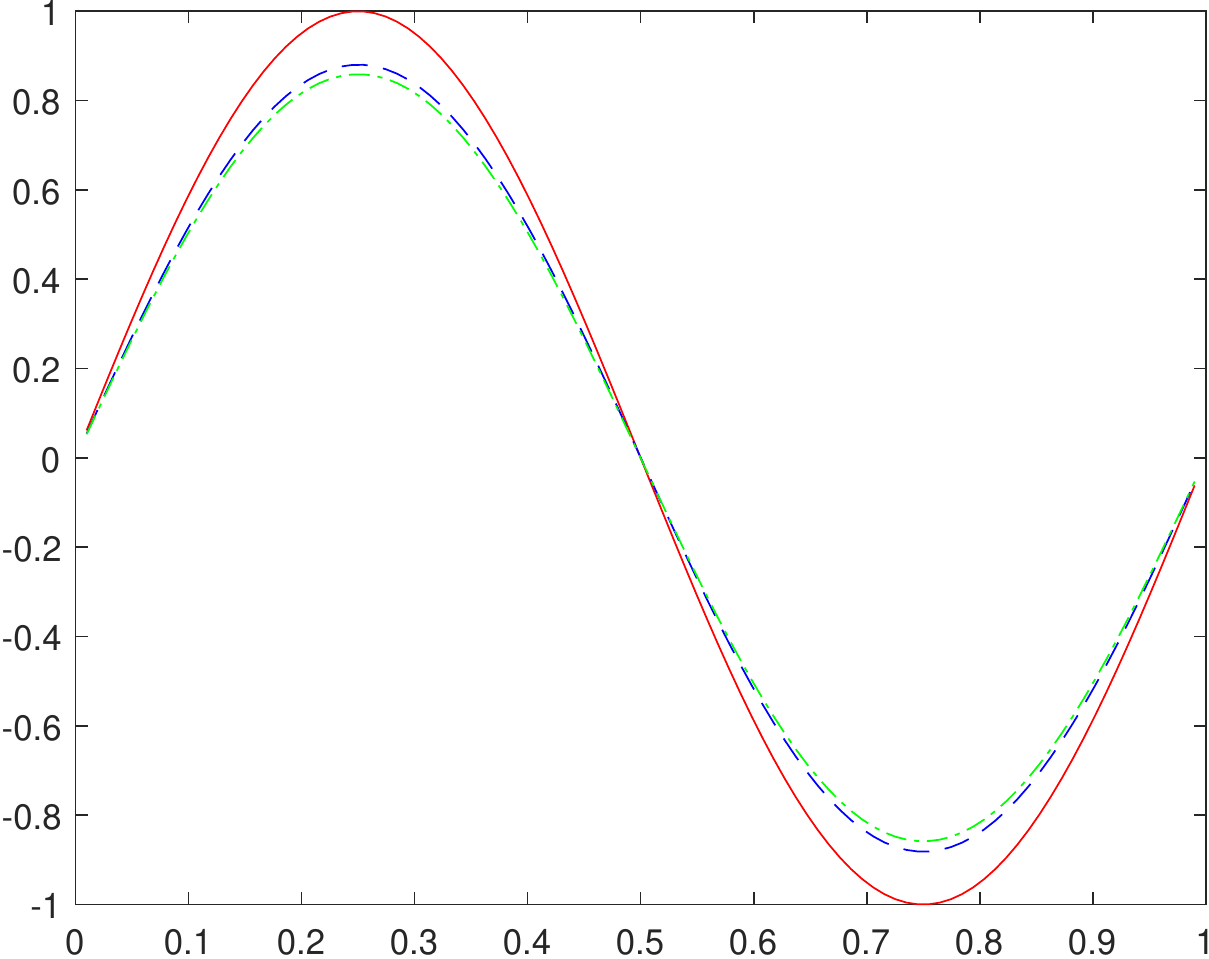}\qquad
\includegraphics[width=3.7cm]{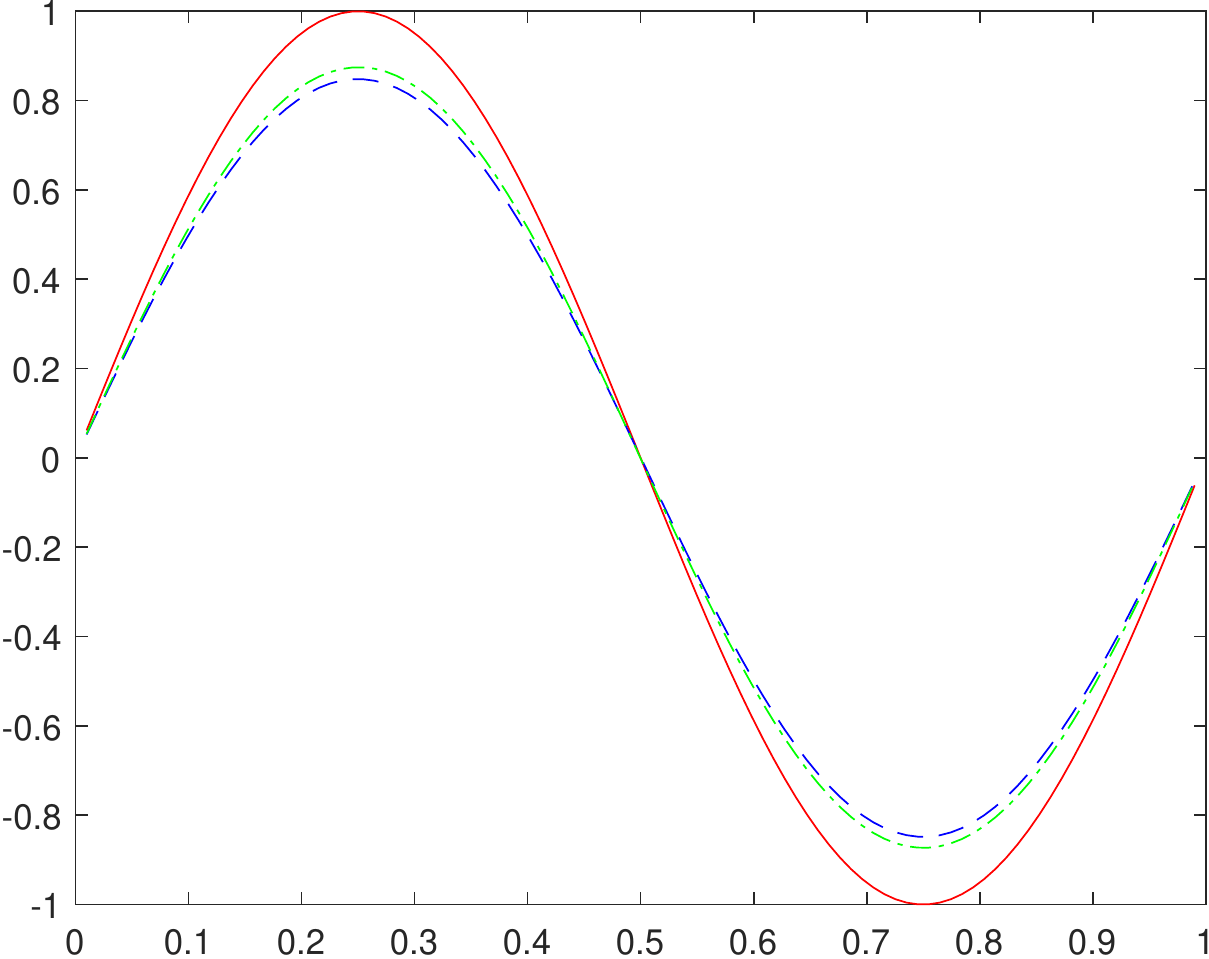}
\caption{\color{blue}{Exact parameter $\theta^\dagger$ (solid red), reconstruction by Landweber-Kaczmarz (dashed blue) and reconstruction by Landweber (dashed dotted green). Left: Reduced, right: AAO setting.}
\label{fig-3}}\vspace{-5pt}
\end{figure}

\begin{figure}[!htb] 
\setlength{\belowcaptionskip}{-5pt}
\centering
\includegraphics[width=3.7cm]{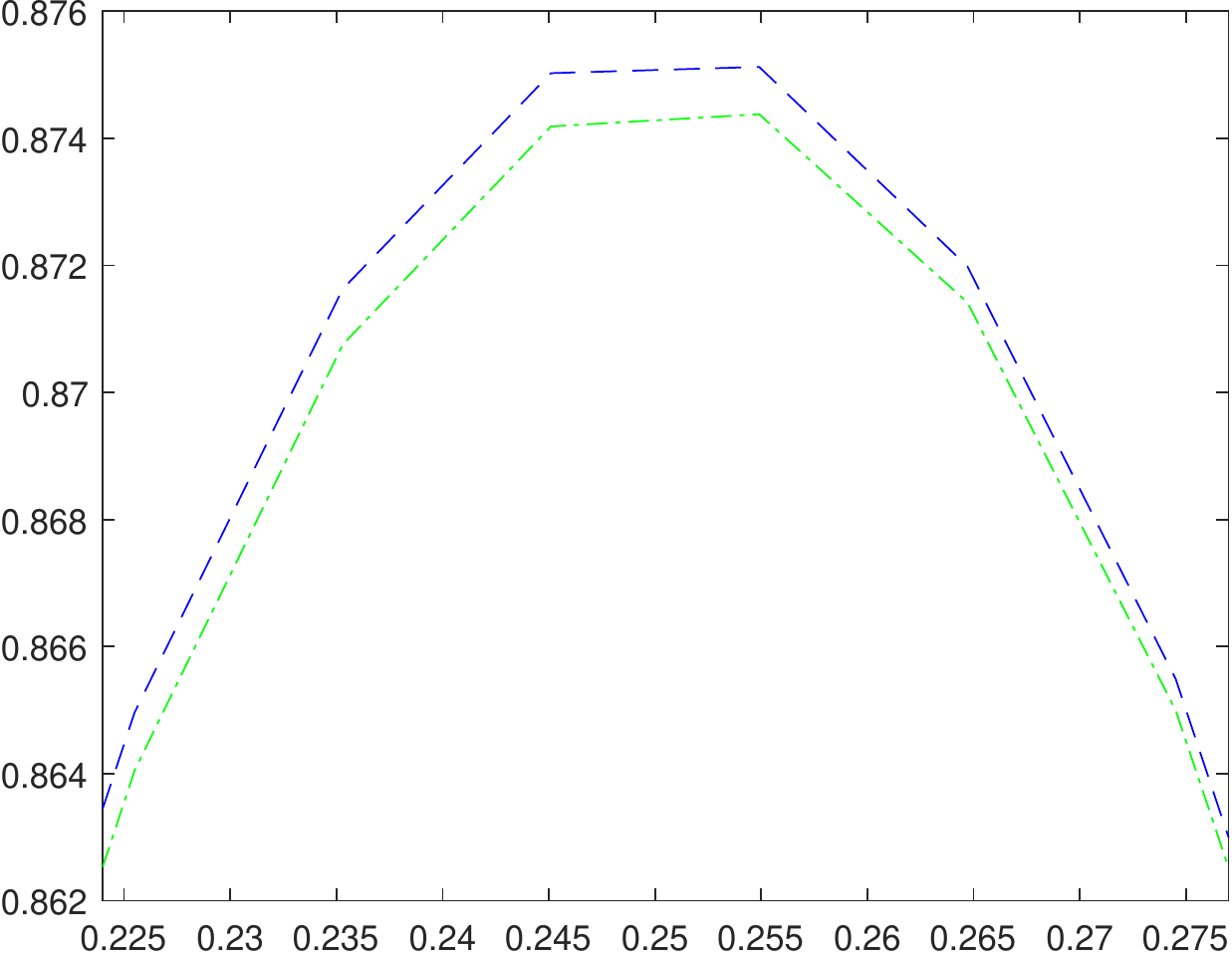}\qquad
\includegraphics[width=3.7cm]{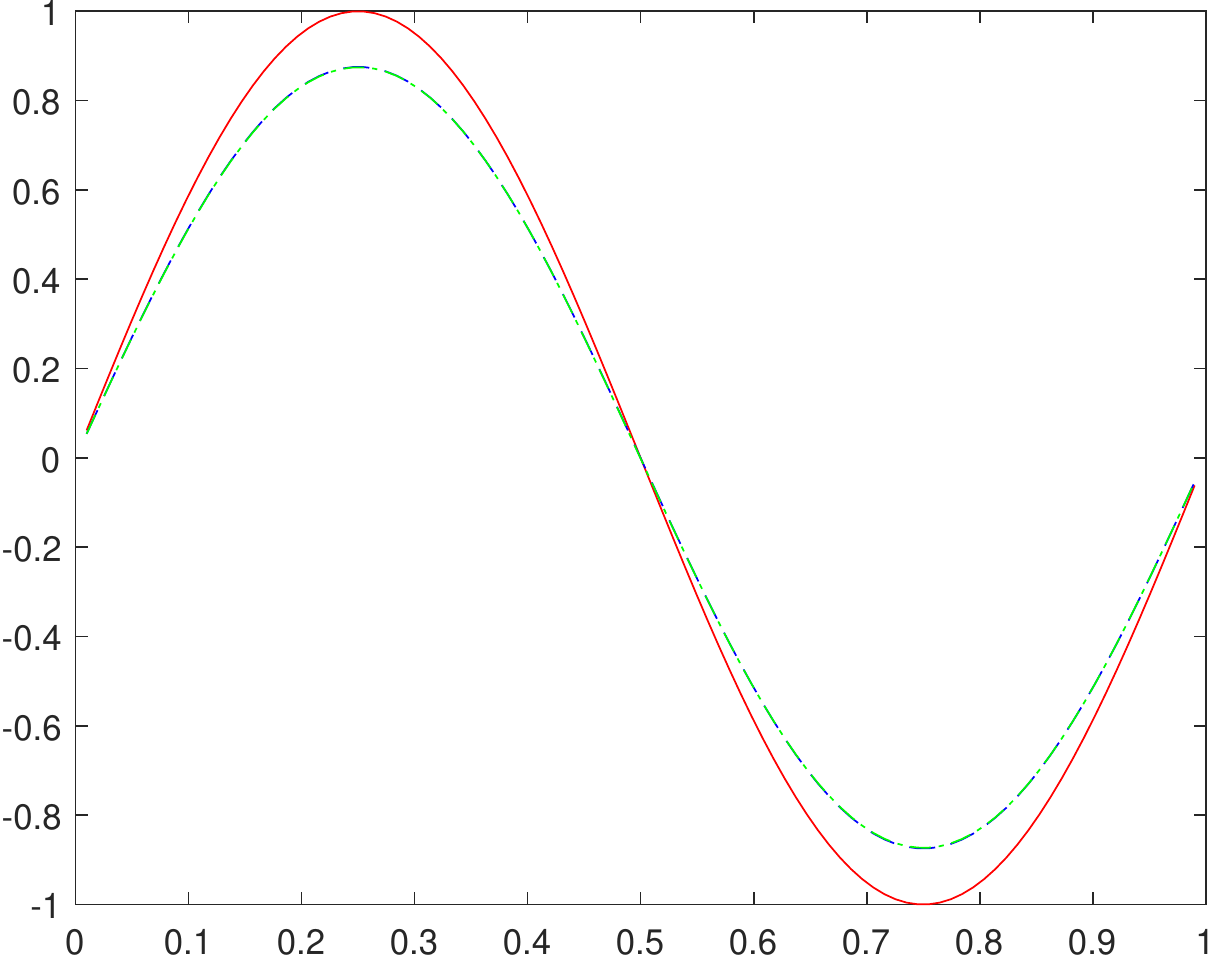}
\caption{\color{blue} Exact parameter $\theta^\dagger$ (solid red), reconstruction by Landweber-Kaczmarz (dashed blue) and reconstruction by Landweber (dashed dotted green). Left: zoom of right, right: AAO setting.}
\label{fig-4}\vspace{-5pt}
\end{figure}

\color{blue}
We turn to a comparison between the Landweber method on the full time horizon and the proposed Landweber-Kaczmarz method. Figure \ref{fig-3} shows that in the Reduced setting the reconstruction by the Landweber-Kaczmarz method is more accurate, but the performance seems to be reversed in the All-at-once setting. The reason for this is that, in the Landweber-Kaczmarz method the initial data appears only in $F_0$, while the Landweber method, as working with a fixed operator $F$ on the whole time line $[0,T]$, employs the initial condition always. This leads to the next test in the All-at-once setting displayed in Figure \ref{fig-4}. In this test, we design the collection of operators $\{F_j\}_{j=0\ldots n-1}$ for Landweber-Kaczmarz method in a slightly different way with the formulation in Remark \ref{AAO-Kaczmarz}, namely, we include the initial data in all $F_j, j=0\ldots n-1$. Figure \ref{fig-4} indicates an enhancement in the computed output from the Landweber-Kaczmarz method, and it now performs slightly better than the Landweber method. Besides, we see the All-at-once result in Figure \ref{fig-4} is more similar to the Reduced result than that in the original Landweber-Kaczmarz method (Figure \ref{fig-1}). This confirms the advantage of the presence of initial data in the All-at-once formulation of the Landweber-Kaczmarz method in this example.

From this, we observe the following
\begin{itemize}
\item Landweber-Kaczmarz performs slightly better than Landweber in both settings.
\item The All-at-once result is more comparable to the Reduced one if the initial data is present in all $F_j$ of the All-at-once formulation.
\item Relying on the idea of time line segmenting, there would be flexible strategies to generate the collection of $F_j$ for Landweber-Kaczmarz method, e.g.:
\begin{itemize}
\item[$\circ$] Method suggested in Remark \ref{AAO-Kaczmarz}
\item[$\circ$] Method suggested in Remark \ref{AAO-Kaczmarz} with $F_j$ defined on $\{0\}\cup[\tau_j, \tau_{j+1}], j=0\ldots n-1$
\item[$\circ$] Evolution of time subinterval, i.e., $F_j$ is defined on $[0, \tau_{j+1}], j=0\ldots n-1$
\end{itemize}
and those strategies affect the runtime accordingly.
\end{itemize}
\color{black}

In all tests, we used the finite difference method to compute the exact state. More precisely, an implicit Euler scheme was employed for time discretization and a central difference quotient was used for space discretization. The numerical integration ran with the trapezoid rule. Gaussian noise was generated to be added to the exact data $y$.

\color{blue}
We also point out that in case of having sufficiently smooth data $y^\delta=u^\delta\in L^2(0,T;H_0^1(\Omega))$, we are able to recover the parameter directly from 
\[ \theta^\delta=\left(\dot{u}^\delta-\Delta u^\delta +{u^\delta}^3\right)(t)\]
if we compute the terms $\dot{u}^\delta$ and $\Delta u^\delta$, e.g., by filtering (cf. Remark \ref{uniqueness}). Getting the output from those in $L^2(0,T;L^2(\Omega))$ and the fact ${u^\delta}^3\in L^2(0,T;L^2(\Omega))$ since $H_0^1(\Omega)\hookrightarrow L^6(\Omega)$, it yields the output $\theta^\delta\in L^2(\Omega)$. Nevertheless, if applying the method developed in the paper, one can deal with much higher nonlinearity, namely, $\Phi(u)=u^5$ as claimed in Remark \ref{smoothdata}. Moreover, as apposed to this direct inversion, our approach also works for partial (e.g. boundary) or time discrete observations.

The source identification problem \eqref{ex1}-\eqref{ex4} examined here can be extended into time-dependent $\theta\in L^2(0,T;X)$ and implemented in the All-at-once as well as the Reduced setting using the Landweber method as discussed in Remark \ref{timedepend}. 
\color{black}
\section{Conclusion and outlook} \label{secOutlook}
In this study, we consider a general evolution system over a finite time line and investigate parameter identification in it by using Landweber-Kaczmarz regularization. We formulate this problem in two different settings: An All-at-once version and a Reduced version. In each version, both cases of full and discrete observations are taken into account. The main ingredients for the regularization method are: differentiability and adjoint of the derivative of the forward operators. Differentiability was proved by mainly basing on the choice of appropriate function spaces and a local Lipschitz continuity condition. Segmenting the time line into several subintervals gives the idea to the application of a Kaczmarz method. A loping strategy is incorporated into the method forming the loping Landweber-Kaczmarz iteration. The shown example proves that the method is efficient in the practically relevant situation of high nonlinearity.

Several questions arise for future research:

We plan to extend the theory to time-dependent parameters. For this purpose, we need to build an appropriate function space for $\theta$ which, for instance, could allow the local Lipschitz continuity condition. In addition, the assumptions for well-posedness of the parameter-to-state $S$ map need to be carefully considered.

Concerning the model, we intend to also study second order in time equations modeling wave phenomena. Rewriting them as first order in time system by introducing another state $\tilde{u}=\dot{u}$, in principle, allows us to use the present formulation. However, an appropriate function space setting for wave type equations requires different tools for showing, e.g., well-definedness of the parameter-to-state map.

In our test problem (\ref{ex1})-(\ref{ex4}), we consider full space observations in order to establish the tangential cone condition. Practically, relevant partial or boundary observations are yet to be tested numerically.

\color{blue}
Regarding numerical implementation of other iterative regularization methods, the difficulty of Newton type methods, which are supposed to give rapid convergence, is the requirement of solving a linear system per iteration step, while this is avoided in the Landweber-Kaczmarz method by the direct use of the Hilbert space adjoints. In the context of this paper, only numerical experiments for Landweber and Landweber-Kaczmarz methods are provided. Numerical implementation and computational tests for other iterative methods will be a subject of future work.

Considering the All-at-once formulation on a infinite time horizon, by setting $\dot{\theta}=0$ and $\bU:=(\theta,u)$, the problem can be written as a dynamical system $\dot{\bU}(t) = (0, f)(t) =: \bF(., \bU)(t), t >0, \bU(0) =(\theta^*, u_0)$, where the exact parameter $\theta^*$, being a time constant function, is supposed to be estimated simultaneously to the time evolution of the system and the data collecting process $y(t)=g(.,u,\theta)(t)=:\bG(\bU)(t)$. This appears to be a link to online parameter identification methods (see e.g. \cite{Kugler1, Kugler}).
The relation between the proposed All-at-once formulation and online parameter identification for time-independent parameters as well as their analysis (possibly by means of a Lopatinskii condition) will be subject of future research.
\color{black}
\section*{Appendix}

\subsection{Regularity result in time-dependent case} \label{Roubicek}
We refer the reader to Regularity Theorem 8.16 in \cite{Roubicek}, which we are using with exactly the same notations. All the equations referred to \enquote{((.))} indicate the ones in the book \cite{Roubicek}.  

\begin{remark}
Observations on the Theorem 8.16 in \cite{Roubicek} 
\begin{itemize}
\item
This proof still holds for the case $A_2$ is time-dependent. The condition ((8.59d)) on $A_2$ could stay fixed or be weakened to $\|A_2(t,v)\|_H \leq C(\gamma(t)+\|v\|_V^{q/2}), \gamma \in L^2(0,T)$, then
\begin{align}
\|u(t)\|_V \leq N\left( 2C\|\gamma\|_{L^2(0,T)} + 2C\sqrt{\frac{c_1}{c_0}}\|u_0\|_H + \|f\|_{L^2(0,T;H)} + \sqrt{|\Phi(u_0)|} \right).
\end{align}
\item
We can slightly relax the constraint on $c_1$ by applying Cauchy's inequality with $\epsilon$ for the first estimate in the original proof. In this way, we get $2c_1T<1$ which can be traded off by the scaling $\frac{1}{\epsilon}$ on the right hand side of ((8.63)).\\ 
If $\Phi: V \rightarrow \R^+$, the assumption on the smallness of $c_1$ can be omitted.

\end{itemize}
\end{remark}

\proof The first expression in the book shows that $A_1$ does not depend on $t$. With the hope of generalizing to time-dependent case,
 our strategy is as follows. \\
First we set
\begin{align*}
& A_1(t,v)=\Phi'_v(t,v), \quad \Phi: [0,T)\times V \rightarrow \R \qquad\qquad\qquad\qquad\\
& \Phi'_v \text{ induces a Nemytskii operator, namely, } \Psi'_u\\
&\text{i.e.}, \Psi'_u(t,u)(t)=\Phi'_v(t,u(t)),
\end{align*}
thus
\begin{align*}
\dupair{A_1(t,u)(t),\frac{d}{dt}u(t)}
&=\dupair{\Phi'_v(t,u(t)),\frac{d}{dt}u(t)} \qquad\qquad\qquad\qquad\\
&=\dupair{\Psi'_u(t,u)(t),\frac{d}{dt}u(t)}\\
&=\Psi'_t(t,u)(t)-\Phi'_t(t,u(t)).
\end{align*}
Looking at the first estimate in the book, we can think of treating $\Psi'_t(t,u)(t)$ as $\Phi'_t(t,u)(t)$ on the left hand side and leaving $\Phi'_t(t,u(t))$ to the right hand side.
 Choosing
\begin{align*}
&\Phi(t,v) \geq c_0\|v\|_V^q - c_1\|v\|_H^2, \qquad\hphantom{l} \forall t \in (0,T) \qquad\qquad\qquad\qquad\\
&\|\Phi'_t(t,v)\|_H \leq \tilde{C}(\tilde{\gamma}(t)+\|v\|_V^q), \quad \tilde{\gamma} \in L^1(0,T)
\end{align*} 
lets us estimate analogously to ((8.62)) and obtain
\begin{align} \label{regularity}
&\|u(t)\|_V \\
&\leq N\left( 2C\|\gamma\|_{L^2(0,T)} + 2\tilde{C}\|\tilde{\gamma}\|_{L^1(0,T)}^\frac{1}{2} + 2(C+\tilde{C})\sqrt{\frac{c_1}{c_0}}\|u_0\|_H + \|f\|_{L^2(0,T;H)} + \sqrt{|\Phi(0,u_0)|} \right) \nonumber
\end{align}
for all $t \in (0,T).$ This completes the proof. \qed

\subsection{Tangential cone condition}\label{TCC}
\begin{itemize}
\item All-at-once setting\medbreak
Since the observation is linear in our example, the condition (\ref{AAO-TCC}) is fulfilled provided that, for every $u,\tilde{u}\in \calB_{2\rho(u_0)}$,
\begin{equation} \label{AAO-TCC1}
\|f(.,\tilde{u},\tilde{\theta}) - f(.,u,\theta) - f'_u(.,u,\theta)(\tilde{u}-u) - f'_\theta(.,u,\theta)(\tilde{\theta}-\theta) \|_{\calW}
\leq c_{tc} \|\tilde{u} - u  \|_{\calY} 
\end{equation}
or, for every $t\in (0,T)$,
\begin{align} \label{AAO-TCC2}
&\|\Phi(\tilde{u})-\Phi(u)-\Phi'(u)(\tilde{u}-u) \|_{V^*}
 \leq c_{tc} \|\tilde{u} - u \|_{Z}. \qquad
\end{align}
Developing the left hand side (LHS) of (\ref{AAO-TCC2}), we have
\begin{align*}
\text{ LHS} 
&= \left\|\int_0^1\int_0^1 \Phi''(u+\sigma\lambda(\tilde{u}-u))(\tilde{u}-u)^2d\lambda d\sigma \right\|_{V^*} \\
&\leq \sup_{0\leq \sigma,\lambda\leq 1}C \left( \int_\Omega\left( \Phi''(u+\sigma\lambda(\tilde{u}-u))(\tilde{u}-u)^2 \right)^\frac{\bar{p}}{\bar{p}-1} dx \right)^\frac{\bar{p}-1}{\bar{p}}\\
&\leq C \|\tilde{u}-u\|_{L^2(\Omega)} \sup_{0\leq \sigma,\lambda\leq 1} \left( \int_\Omega\left( |u+\sigma\lambda(\tilde{u}-u)|^{\gamma-2}(\tilde{u}-u) \right)^\frac{2\bar{p}}{\bar{p}-2} dx \right)^\frac{\bar{p}-2}{2\bar{p}}\\
&\leq C \|\tilde{u}-u\|_{L^2(\Omega)}\|\tilde{u}-u\|_V \left( \|\tilde{u}^{\gamma-2}\|_{L^{\frac{2\bar{p}}{\bar{p}-4}}}+\|u^{\gamma-2}\|_{L^{\frac{2\bar{p}}{\bar{p}-4}}} \right)\\
&\leq C \|\tilde{u}-u\|_{L^2(\Omega)}\|\tilde{u}-u\|_\calU \left( \|\tilde{u}\|_\calU+\|u\|_\calU \right) \quad\qquad\qquad\qquad\qquad\\
&\leq c_{tc} \|\tilde{u}-u\|_{L^2(\Omega)},
\end{align*}
\textcolor{blue}{where the generic constant $C$ may take different values whenever it appears.} The tangential cone coefficient $c_{tc}$, which depends only on $c_{H^1\rightarrow L^{\bar{p}}}, c_{PF}, \gamma$ and $T$, is sufficiently small if $u$ is sufficiently close to $\tilde{u}$ and $\gamma\leq \frac{\bar{p}}{2}$.

\item Reduced setting\medbreak
We need to verify that, for all $\theta, \tilde{\theta}\in\calB_{2\tilde{\rho}(\theta_0)}$,
\begin{align} \label{Re-TCC1}
\|S(\tilde{\theta})-S(\theta)-v\|_\calY\leq \tilde{c}_{tc}\|S(\tilde{\theta})-S(\theta)\|_\calY,
\end{align}
where $v$ solves (\ref{Re-TTCv}). Letting $\xi=\tilde{\theta}-\theta$ then $v=u^\xi$ solves the sensitivity equation (\ref{S'eq}), and by denoting $\tilde{u}:=S(\tilde{\theta}), u:=S(\theta)$, (\ref{Re-TCC1}) becomes
\begin{align*}
\|\tilde{u}-u-u^\xi\|_\calY\leq \tilde{c}_{tc}\|\tilde{u}-u\|_\calY. \quad
\end{align*}
Denoting $\tilde{u}-u-u^\xi=:\tilde{v}_1=\tilde{v}_{\epsilon=1}$, with $\tilde{v}_\epsilon$ as in (\ref{Re-vtilde}), and using the fact that $\tilde{v}_1(0)=0$, we deduce
\begin{align*}
\|\tilde{v}_1\|_{L^2(0,T;L^2(\Omega))} &\leq C\|\tilde{v}_1\|_{W^{1,2,2}(0,T;V,V^*)}\leq CC^\theta\|r_1\|_{L^2(0,T)}\qquad\qquad\qquad\\
&\leq CC^\theta c_{tc}\|\tilde{u}-u\|_{L^2(0,T;L^2(\Omega))} \\
&=:\tilde{c}_{tc}\|\tilde{u}-u\|_{L^2(0,T;L^2(\Omega))},
\end{align*}
where $r_1=r_{\epsilon=1}$, and its $L^2(0,T)$-norm is the left hand side of $(\ref{AAO-TCC1})$.
\end{itemize}

\begin{remark} \label{smoothdata}
In three dimensions, $\gamma=5$ is achievable with $\calY=L^2(0,T;V)$, however the realistic data space is $\calY=L^2(0,T;L^2(\Omega))$.
\end{remark}

\section*{Acknowledgment}
The author wishes to thank \textcolor{blue}{her supervisor} Barbara Kaltenbacher, Alpen-Adria-Universit\" at Klagenfurt, for fruitful and inspiring discussion.\\
\textcolor{blue}{The author would like to thank the reviewers for the insightful comments leading to an improvement of the manuscript.}\\
The author gratefully acknowledges partial support by the Karl Popper Kolleg ``Modeling-Simulation-Optimization'', funded by the Alpen-Adria-Universit\" at Klagenfurt and by the Carinthian Economic Promotion Fund (KWF).

\section*{References}
\bibliographystyle{siam}
\nocite{*}
\bibliography{Reference-OLD}

\end{document}